\documentclass[12pt]{amsart}
\usepackage{amssymb}
\usepackage{amsmath}
\usepackage{versions}
\usepackage{color}
\numberwithin{equation}{section}
\newtheorem{theorem}{Theorem}[section]
\newtheorem{lemma}[theorem]{Lemma}
\newtheorem{corollary}[theorem]{Corollary}
\newtheorem{example}[theorem]{Example}
\newtheorem{examples}[theorem]{Examples}
\newtheorem{definition}[theorem]{Definition}
\newtheorem{proposition}[theorem]{Proposition}
\newtheorem{conjecture}[theorem]{Conjecture}
\newtheorem{observation}[theorem]{Observation}
\newtheorem{remark}[theorem]{Remark}

\newtheorem{question}[theorem]{Question}
\newcommand\beq{\begin{equation}}
\newcommand\eeq{\end{equation}}
\newcommand\re{\mathrm {Re~}}

\newcommand\ii{\mathrm i}
\newcommand\al{\alpha}

\newcommand\de{\delta}
\newcommand\la{\lambda}
\newcommand\up{\upsilon}
\newcommand\ups{\upsilon}
\newcommand{\ph}{\varphi}

\newcommand\ta{\theta}
\newcommand\D{\mathbb D}
\newcommand\T{\mathbb T}
\newcommand\C{\mathbb C}
\newcommand\G{\mathbb{G}}
\newcommand\R{\mathbb R}
\newcommand\e{\mathrm e}
\newcommand\N{\mathbb N}
\newcommand\E{\mathcal E}
\newcommand\B{\mathcal B}
\newcommand\V{\mathcal V}

\newcommand\id{\rm id}
\newcommand\diag{\mathrm {diag~}}
\newcommand\df{\stackrel{\rm def}{=}}
\newcommand\Aut{\mathrm{Aut}~}
\newcommand\Schur{\mathcal{S}}

\newcommand\nn{\nonumber}

\newcommand\black{\color{black}}

\DeclareMathOperator\hol{Hol}

\def\qed{\hfill $\square$ \vspace{3mm}}
\let\phi\varphi
\numberwithin{equation}{section}
\begin{document}

\title{Extremal holomorphic maps and the symmetrised bidisc}
\author{Jim Agler, Zinaida A. Lykova and N. J. Young}
\date{30th July, 2011}

\begin{abstract} 
 We introduce  the class of $n$-extremal holomorphic maps, a class that generalises both finite Blaschke products and complex geodesics, and apply the notion to the finite interpolation problem for analytic functions from the open unit disc into the symmetrised bidisc $\Gamma$.  We show that a well-known necessary condition for the solvability of such an interpolation problem is not sufficient whenever the number of interpolation nodes is $3$ or greater.
We introduce a sequence $\mathcal{C}_\nu, \nu \geq 0,$ of necessary conditions for solvability, prove that they are of strictly increasing strength and show that $\mathcal{C}_{n-3}$ is insufficient for the solvability of an $n$-point problem for $n\geq 3$.  We propose the conjecture that condition $\mathcal{C}_{n-2}$ is necessary and sufficient for the solvability of an $n$-point interpolation problem for $\Gamma$ and we explore the implications of this conjecture.

We introduce a classification of rational $\Gamma$-inner functions, that is, analytic functions from the disc into $\Gamma$ whose radial limits at almost all points on the unit circle lie in the distinguished boundary of $\Gamma$.  The classes are related to $n$-extremality and the conditions $\mathcal{C}_\nu$; we prove numerous strict inclusions between the classes.
\end{abstract}


\thanks{The first author was partially supported by National Science Foundation Grant on  Extending Hilbert Space Operators DMS 1068830. The third author was partially supported by the UK Engineering and Physical Sciences Research Council grant EP/J004545/1.}
\maketitle

\section*{Contents} \label{contents}

\ref{newintro}. Introduction \hfill  Page \pageref{newintro}

\ref{n_extremal}. $n$-extremal holomorphic maps   \hfill \pageref{n_extremal}

\hspace*{1cm}\ref{geos2extrem}. Complex geodesics and $2$-extremals \hfill \pageref{geos2extrem}

\hspace*{1cm}\ref  {universal}.  Universal Carath\'eodory sets and $n$-extremal maps \hfill \pageref{universal}

\ref{symmetrised}.  The symmetrised bidisc $\Gamma$\hfill \pageref{symmetrised}

\ref{interp}. Interpolation in $\hol(\D,\Gamma)$ and a conjecture \hfill\pageref{interp}

\ref{extremality}. Extremality in Condition $\mathcal{C}_\nu$ \hfill\pageref{extremality}

\ref{Gamin}. $\Gamma$-inner functions \hfill \pageref{Gamin}

\hspace*{1cm} \ref{newOld}.   New $\Gamma$-inner functions from old \hfill\pageref{newOld}

\hspace*{1cm} \ref{RatGamIn}. Rational $\Gamma$-inner functions \hfill \pageref{RatGamIn}

\ref{Eclasses}. The classes $\E_{\nu k}$ \hfill\pageref{Eclasses}

\hspace*{1cm}\ref{phasar}.  Phasar derivatives \hfill\pageref{phasar}

\hspace*{1cm}\ref{cancels}.  Cancellations and the classes $\E_{\nu k}$\hfill\pageref{cancels}

\ref{superficial}. Superficial $\Gamma$-inner functions  and the classes $\E_{\nu 1}$\hfill\pageref{superficial}

\ref{kExtremEnuk}. The classes $\E_{\nu k}$ and $k$-extremals \hfill\pageref{kExtremEnuk}

\ref{geodesics}. Complex geodesics of $\G $ and the classes $\E_{\nu 2}$\hfill\pageref{geodesics}

\ref{connections-C-E}. Condition $\mathcal{C}_\nu$ and the classes $\E_{\nu k}$\hfill\pageref{connections-C-E}

\ref{inequations-for-E}. Inequations for the classes $\E_{\nu k}$\hfill\pageref{inequations-for-E}

\ref{summary-E-cl}. Table of relations between the classes $\E_{\nu k}$\hfill\pageref{summary-E-cl}

\ref{conclud}. Concluding reflections \hfill\pageref{conclud}

 References \hfill\pageref{bibliog}

\section{Introduction} \label{newintro}

One result in G. Pick's seminal paper \cite{Pick} of 1916  states that finite Blaschke products of prescribed degree are characterised by a certain extremal property.  Let us say that  an analytic self-map $h$ of the open unit disc $\D$ is {\em $n$-extremal} if 
for every $n$-point subset $\Lambda$ of $\D$, there do not exist an  $r>1$ and a holomorphic function $f: r\D\to \D$ such that $f$ and $h$ agree on $\Lambda$.  A version of Pick's result can be formulated: the $n$-extremal holomorphic self-maps of $\D$ are precisely the Blaschke products of degree at most $n-1$.

A similar notion of extremality, but with $n$ equal to $2$, occurs in the theory of hyperbolic complex spaces introduced by S. Kobayashi and described in his book \cite{Ko98}.  In this context one studies the geometry and function theory of a domain $\Omega\subset \C^d$ with the aid of $2$-extremal holomorphic maps from $\D$ to $\Omega$.  The notion of $n$-extremal map makes sense, however, in much greater generality. We consider it for elements of the space $\hol(\Omega, E)$ of holomorphic maps from a domain $\Omega$ to a subset $E$ of $\C^N$ (Definition \ref{extSolve}).  

A prominent theme in hyperbolic complex geometry is a kind of duality between $\hol(\D,\Omega)$ and $\hol(\Omega,\D)$, typified by the celebrated theorem of L. Lempert \cite{Le86},  which in our terminology asserts that if $\Omega$ is convex then every $2$-extremal map belonging to $\hol(\D,\Omega)$ is a complex geodesic of $\Omega$ (that is, has an analytic left inverse).  This duality plays a role in the present paper too.

Since $n$-extremal maps simultaneously generalise both Blaschke products of prescribed degree and complex geodesics, they surely constitute a significant class.  We have encountered them in attempting to solve a certain interpolation problem, a special case of the {\em $\mu$-synthesis problem}, which arises in control engineering.  We make some remarks about this application in Section \ref{conclud} below; see also \cite{BFT1,NJY11} for more on this topic.  The problem led us to investigate \cite{AY1,AY01} a special domain in $\C^2$, the {\em symmetrised bidisc}, defined to be the set
\[
\G \df \{(z+w,zw): z,w\in\D\}.
\]
The rich and surprising geometry of this domain and its higher-dimensional analogues has subsequently been elaborated by many authors  (for example \cite{costara,EZ,JP,NiPfZw,ogle}).

In this paper we focus on the finite interpolation problem for $\hol(\D,\Gamma)$, where $\Gamma$ is the closure of $\G$: given $n$ interpolation nodes in $\D$ and $n$ target points in $\Gamma$, we wish to determine whether there exists a function in $\hol(\D,\Gamma)$ that satisfies the corresponding interpolation conditions.  If $\Gamma$ were replaced by the closure of a Cartan domain then the beautiful and far-reaching classical Nevanlinna-Pick theory would apply, but the domain $\G$ is inhomogeneous and at present there is no satisfactory criterion for the solvability of interpolation problems in $\hol(\D,\Gamma)$.  There {\em is} a well-known necessary condition, which we call $\mathcal{C}_0$, and which is numerically practicable for modest $n$  (one must check for positivity a one-parameter family of Hermitian $n\times n$ matrices).  Condition $\mathcal{C}_0$ is also sufficient for solvability when $n=2$; one of our principal results is that condition $\mathcal{C}_0$ is {\em not} sufficient for  solvability if $n\geq3$.  Accordingly, we introduce a sequence of necessary conditions $\mathcal{C}_\nu$, for $\nu=0,1,2,\dots$ (Definition \ref{C_nu}).  We prove that this sequence of conditions is of strictly increasing strength, from which it follows that none of the $\mathcal{C}_\nu$ is sufficient for all finite interpolation problems in $\hol(\D,\Gamma)$.  Nevertheless, it seems possible that $\mathcal{C}_{n-2}$ is sufficient for the  solvability of an $n$-point problem: we conjecture that this is indeed so (the {\em $\Gamma$-interpolation Conjecture}, Section \ref{interp}), and we explore the consequences of this conjecture.

To obtain our results we analyse {\em $\Gamma$-inner functions}: these are maps $h\in\hol(\D,\Gamma)$ whose radial limits almost everywhere on the unit circle $\T$ lie in the distinguished boundary of $\Gamma$.
A good understanding of rational $\Gamma$-inner functions is likely to play a part in any future solution of the finite interpolation problem for $\hol(\D, \Gamma)$, since such a problem has a solution if and only if it has a rational $\Gamma$-inner solution  (for example, \cite[Theorem 4.2]{Cost05}). 
We introduce an array $\E_{\nu n}$ of classes of rational $\Gamma$-inner functions that are closely related both to $n$-extremal maps and to the conditions $\mathcal{C}_\nu$.
We say that a  function $f =(s,p) $ is in  $\E_{\nu n}$ if $f\in\hol(\D,\Gamma)$ is rational and there exists $m \in \B l_{\nu}$ such that 
$$
\frac{2 m p -s}{2- m s} \in \B l_{n-1}.
$$
Here $\B l_n$ is the set of Blaschke products of degree at most $n$.  We show that any function in $\E_{\nu n}$ either maps into the topological boundary of $\Gamma$ or is $n$-extremal, while if the $\Gamma$-interpolation Conjecture holds, then any rational $n$-extremal $\Gamma$-inner function belongs to $\E_{n-2,n}$.  We obtain numerous strict inclusions between $\E$ classes, which are summarised in a table in Section \ref{summary-E-cl}.

Here is some terminology and notation.

The closed unit disc $\{z : |z| \le 1 \}$ will be denoted by $\Delta$.  The closure of a set $S$ in a topological space will be written $S^-$.
We denote by $\mathbb{T}$ the unit circle,  by $H^2$ the Hardy Hilbert space on $\mathbb{D}$ and by $K$ the Szeg\H{o} kernel:
$$
K_\lambda (z) = K(z, \lambda) = (1-\overline \lambda z)^{-1},
\quad \lambda, z\in \mathbb{D}.
$$
For any domain $\Omega$ and any set $E \subset \C^{N}$, we denote by ${\rm Hol}(\Omega, E)$ the set of analytic functions from $\Omega$ into $E$.
The {\em Schur class} $\Schur$ is the class ${\rm Hol}(\D,\Delta)$ of functions analytic and
bounded by 1 in $\mathbb{D}$.  For a function $f$ on a subset of the complex plane $\mathbb{C}$
we write
$$
\bar f(z) = \overline{f(z)},\qquad f^\vee(z) = (f(\bar z))^-.
$$

For $\alpha \in \C$ we write
$$
B_\alpha(z) = \frac{z-\alpha}{1-\overline \alpha z}.
$$
In the event that $\al\in\D$ the rational function $B_\al$ is called
a {\em Blaschke factor}.
A {\em M\"obius function} is a function of the form $cB_\alpha$
for some $\alpha \in \mathbb{D}$ and $c\in \mathbb{T}$.
The set of all M\"obius functions is the automorphism group 
$\Aut \mathbb{D}$ of $\mathbb{D}$.  

We denote by $d(f)$ the degree of a rational function $f$ of one variable -- that is, the maximum of the degrees of the numerator and denominator  in an expression of $f$ as a ratio of coprime polynomials.  

\section{$n$-extremal holomorphic maps} \label{n_extremal}

  Roughly speaking, a holomorphic map $h$ between domains is $n$-extremal if its restriction to any $n$-point set yields interpolation data that are solvable, but only just.  More precisely:

\begin{definition}\label{extSolve} 
Let $\Omega$ be a domain, let $E\subset \C^N$, let $n\geq 1$, let $\la_1,\dots,\la_n$ be distinct points in $\Omega$ and let $z_1,\dots,z_n \in E$.  We say that the interpolation data 
\[
\la_j \mapsto z_j : \Omega \to E,   \quad j=1,\dots, n,
\]
 are {\em extremally solvable} if there exists a map $h\in \hol(\Omega,E)$ such that $h(\la_j)=z_j$ for $j=1,\dots,n$, but, for any open neighbourhood $U$ of the closure of $\Omega$, there is no $f\in\hol(U,E)$ such that $f(\la_j)=z_j$ for $j=1,\dots,n$.

We say further that $h\in \hol(\Omega,E)$ is {\em $n$-extremal (for $\hol(\Omega,E)$)} if, for all choices of $n$ distinct points $\la_1,\dots,\la_n$ in $\Omega$, the interpolation data
\[
\la_j \mapsto h(\la_j): \Omega \to E, \quad j=1,\dots ,n,
\]
are extremally solvable.
\end{definition}
There are no $1$-extremal holomorphic maps, so we shall always suppose that $n\geq 2$.

As we mentioned in the Introduction, Pick showed that a function $f$ is $n$-extremal for the Schur class $\Schur=\hol(\D,\Delta)$ if and only if $f\in \B l_{n-1}$.   Thus $n$-extremals may be regarded as an analogue of the Blaschke products of degree at most $n-1$.

It is evident that the notion of an $n$-extremal holomorphic map applies very generally, but in this paper we shall be mainly concerned with $n$-extremals for $\hol(\D,\Gamma)$.  We shall however point out some simple general properties of $n$-extremals.

Firstly,  if $h$ is an $n$-extremal for $\hol(\Omega,E)$ then it is also an $m$-extremal for all $m\geq n$.  This is immediate from the definition.

Secondly we note a property of compositions of holomorphic maps: if $g\circ \alpha$ is $n$-extremal then so is $\alpha$.
\begin{proposition}\label{compose}
If $E$ is a set in $\C^N$,  $\Omega_1,\, \Omega_2$ are domains and $g\circ \alpha$ is $n$-extremal for some $\alpha\in\hol(\Omega_1,\Omega_2)$ and $g\in\hol(\Omega_2,E)$ then $\alpha$ is $n$-extremal.  
\end{proposition}
\begin{proof}
For if $\alpha\in\hol(\Omega_1,\Omega_2)$ is {\em not} $n$-extremal then there exist distinct points $\la_1,\dots,\la_n\in\Omega_1$, an open neighbourhood $U$ of the closure of $\Omega_1$ and an $f\in\hol(U,\Omega_2)$ such that $f$ and $ \alpha$ agree on $\la_1,\dots,\la_n$.  Then $g\circ f\in \hol(U,E)$ and $g\circ \alpha$ and $g\circ f$ agree on $\la_1,\dots,\la_n$. This shows that $g\circ \alpha$ is not $n$-extremal. 
\end{proof}

Thirdly we consider the question of the holomorphic invariance of the notion of $n$-extremal.  If $g$ is $n$-extremal in $\hol(\Omega,E)$ and $\al$ is an automorphism of $\Omega$, is it the case that $g\circ\al$ is $n$-extremal?  If $\al$ extends holomorphically to a neighbourhood of the closure of $\Omega$ then the answer is yes.
\begin{proposition}\label{holinv}
Let $\alpha:\Omega_1\to\Omega_2$ be a biholomorphic map of domains that extends to a biholomorphic map from an open neighbourhood $U$ of $\Omega_1^-$ to an open neighbourhood $V$ of $\Omega_2^-$.  If $g$ is $n$-extremal in $\hol(\Omega_2,E)$ then $g\circ \alpha$ is $n$-extremal in $\hol(\Omega_1,E)$.
\end{proposition}
\begin{proof}
Suppose that $g\circ \alpha$ is not $n$-extremal: there exist distinct points $\la_1,\dots,\la_n \in \Omega_1$, an open neighbourhood $U_1$ of $\Omega_1^-$ and a function $f\in\hol(U_1,E)$ such that
\[
f(\la_j) = g\circ \alpha(\la_j) \quad \mbox{ for } j=1,\dots,n.
\]
By hypothesis $\alpha$ extends biholomorphically to an open neighbourhood $U$ of $\Omega_1^-$, and therefore maps $U\cap U_1$ onto an open neighbourhood $V_1$ of $\Omega_2^-$.  Thus $f\circ \alpha^{-1} \in\hol(V_1,E)$, and $f\circ \alpha^{-1}$ agrees with $g$ at the $n$ points $\alpha(\la_1), \dots, \alpha(\la_n)$.  Hence $g$ is not $n$-extremal.
\end{proof}

However, for general $\Omega$ and automorphism $\al$ there is no reason to expect $\al$ to extend even continuously to the closure of $\Omega$, and it therefore seems unlikely that $n$-extremality is preserved under composition on the right with automorphisms in general.

There is a dual result to Proposition \ref{holinv}: it shows that $n$-extremality is better behaved with respect to composition on the left by an automorphism. Let us say that $\alpha:E_1 \to E_2$ is an isomorphism, for any pair of sets $E_1$, $E_2$, if $\alpha$ is bijective and $\alpha, \alpha ^{-1}$ are complex-differentiable on $E_1$, $E_2$ respectively.

\begin{proposition}\label{hol_isom} Let $\Omega$ be a domain and let $\alpha:E_1 \to E_2$ be an isomorphism. Then $g \in \hol(\Omega,E_1)$ is $n$-extremal 
if and only if $\alpha \circ g$ is $n$-extremal in $\hol(\Omega,E_2)$.
\end{proposition}
\begin{proof}
Suppose that $g$ is $n$-extremal but  $\alpha \circ g$ is not $n$-extremal.
Then there exist distinct points $\la_1,\dots,\la_n \in \Omega$, an open neighbourhood $U$ of $\Omega^-$ and a function $f\in\hol(U,E_2)$ such that
\[
f(\la_j) = \alpha \circ g (\la_j) \quad \mbox{ for } j=1,\dots,n.
\]
Since $\alpha^{-1}$ is complex-differentiable, $ \alpha^{-1} \circ f \in\hol(U,E_1)$ (the property of complex-differentiability is vacuous at any isolated point of $E_2$, but this does not matter). Furthermore 
\[
 \alpha^{-1} \circ f (\la_j)= g(\la_j)  \quad \mbox{ for } j=1,\dots,n,
\]
contrary to the hypothesis that $g$ is $n$-extremal. Thus $g$  $n$-extremal
implies that $\alpha \circ g$ is  $n$-extremal. On applying this result to 
 $\alpha \circ g \in \hol(\Omega,E_2)$ and $
\alpha^{-1}:E_2 \to E_1$ we obtain the converse statement.
\end{proof}

Our fourth observation relates to the question: for which $n$ and $\Omega$ is the identity map $\id_\Omega$ $n$-extremal?  We do not know any domain for which the identity map is not $2$-extremal, and the following proposition gives a significant class of domains for which $\id_\Omega$ {\em is} $2$-extremal.

\begin{proposition} \label{ident2ext}
Let $\Omega$ be a bounded domain in $\C^N$ such that, for every pair of distinct points $z_1,z_2\in\Omega$ there is a {\em rational} Kobayashi extremal function for $z_1,z_2$.  Then $\id_\Omega$ is $2$-extremal.
\end{proposition}
We defer the explanation and proof of this statement to the next subsection.

\subsection{Complex geodesics and $2$-extremals} \label{geos2extrem}

The $n$-extremal holomorphic maps also generalise a class of functions that are important in complex geometry.  Consider any domain $\Omega\subset \C^N$.   Let $\delta_\Omega$ denote the Lempert function of $\Omega$, defined for $z_1,\ z_2 \in\Omega$ by
\begin{align}\label{defdelta}
\delta_\Omega( z_1,z_2) \df \inf\{\rho(\mu_1,\mu_2):  &\mbox{ there exists } f\in\hol(\D,\Omega) \mbox{ such that }\\ & f(\mu_1)=z_1, f(\mu_2)=z_2\} \nn
\end{align}
where $\rho$ denotes the pseudohyperbolic distance on $\D$,
\[
\rho(z,w) = \left| \frac{z-w}{1- \bar{w} z} \right|.
\]

Let us say that $h\in\hol(\D,\Omega)$ is a {\em Kobayashi disc} in $\Omega$ if, for every $\la_1,\la_2 \in\D$,
\beq\label{kobadisc}
 \delta_\Omega( h(\la_1), h(\la_2)) =\rho(\la_1,\la_2) 
\eeq
(Kobayashi \cite[Chapter 4, Section 6]{Ko98} calls $h$ an extremal disc if equation \eqref{kobadisc} holds for {\em some} pair of distinct points $\la_1,\la_2\in\D$).

\begin{proposition}\label{koba2ext}
A function $h\in\hol(\D,\Omega)$ is $2$-extremal if and only if $h$ is a Kobayashi disc in $\Omega$.
\end{proposition}
\begin{proof}
$\Leftarrow$  Let $h$ be a Kobayashi disc in $\Omega$.  Suppose that $h$ is not $2$-extremal: then there exist distinct $\la_1,\la_2\in\D$, a real number $r>1$ and $f\in\hol(r\D,\Omega)$ such that $f(\la_1)=h(\la_1)$ and $ f(\la_2)=h(\la_2)$.  Define
\[
f_r\in\hol(\D,\Omega): \la \mapsto f(r\la).
\]
We have $f_r(\la_j/r) = f(\la_j) = h(\la_j)$, and so,
since $h$ is a Kobayashi disc,
\[
\rho(\la_1,\la_2) \leq \rho(\la_1/r,\la_2/r).
\]
Hence
\[
|r^2-\bar\la_2\la_1| \leq r|1-\bar\la_2 \la_1|.
\]
On squaring and expanding we find that
\[
r^2(r^2-1) \leq (r^2-1)|\la_1\la_2|^2,
\]
which is a contradiction since $r>1$ and $\la_1,\ \la_2 \in \D$.  Thus $h$ is $2$-extremal.

\noindent $\Rightarrow$  Let $h$ be $2$-extremal.  Suppose that $h$ is not a Kobayashi disc in $\Omega$: then there exist $\la_1,\la_2\in\D$ such that
\[
\delta_\Omega(h(\la_1),h(\la_2)) < \rho(\la_1,\la_2).
\]
Hence there are $\mu_1,\mu_2\in\D$ and $f\in\hol(\D,\Omega)$ such that $f(\mu_j)=h(\la_j), \, j=1,2$ and
\[
\rho(\mu_1,\mu_2) < \rho(\la_1,\la_2).
\]
By composing $f, h$ with automorphisms of $\D$ we can arrange that $\la_1=\mu_1=0$ and
\[
0 < \mu_2 < \la_2 < 1.
\]
Let $r=\la_2/\mu_2 > 1.$  Define 
\[
g\in\hol(r \D,\Omega): \la \mapsto f(\la/r).
\]
We have $g(0)=f(0)=h(0)$ and
\[
g(\la_2)= f(\la_2/r)=f(\mu_2)=h(\la_2).
\]
This contradicts the hypothesis that $h$ be $2$-extremal.  Thus $h$ is a Kobayashi disc in $\Omega$. 
\end{proof}

For a large class of domains, the Kobayashi discs coincide with the complex geodesics.  We recall that
an analytic  function $h: \D \to \Omega $ is called a {\em complex geodesic of $\Omega$} if there exists 
an analytic left inverse $g : \Omega \to \D$ of $h$.

There is a dual notion to the Lempert function: the {\em Carath\'{e}odory pseudodistance}  $C_\Omega$ on $\Omega$ is given by
\[
C_{\Omega}(z_1, z_2) = \sup\{ \rho( F(z_1), F(z_2)): F\in \hol( \Omega, \D)\}.
\]
Any function $F\in\hol(\Omega,\D)$ for which the supremum on the right-hand side is attained is called  a {\em Carath\'eodory extremal function} for $z_1, z_2$ and the domain $\Omega$. 

\begin{corollary}\label{compgeos}
Let $\Omega$ be a domain for which $\delta_\Omega=C_\Omega$.  The $2$-extremals for $\hol(\D,\Omega)$ coincide with the complex geodesics of $\Omega$.
\end{corollary}

For in such domains the Kobayashi discs coincide with the complex geodesics \cite[Corollary 4.6.2]{Ko98}.   

By Lempert's theorem \cite{Le86}, \cite[Theorem 4.8.13]{Ko98}  the equality $\delta_\Omega=C_\Omega$ holds for convex domains, but there are also nonconvex domains for which it holds, including the symmetrised bidisc $\G$, which is not isomorphic to any convex domain \cite{costara}.

Let us return to Proposition \ref{ident2ext}.  A {\em Kobayashi extremal function} for a pair of distinct points $z_1, \ z_2$ in a domain $\Omega$ is a function $h\in\hol(\D,\Omega)$ which is extremal for the infimum in equation \eqref{defdelta}; thus, for such an $h$, there exist $\la_1, \ \la_2 \in\D$ such that 
\beq \label{proph}
h(\la_1)=z_1, \quad h(\la_2)=z_2\quad \mbox{ and } \quad \rho(\la_1,\la_2)= \delta_\Omega(z_1,z_2).
\eeq
{\bf Proof of Proposition \ref{ident2ext}.}    Suppose $\id_\Omega$ is not $2$-extremal.  Then there exist $z_1, z_2 \in\Omega$, an open neighbourhood $U$ of the closure of $\Omega$ and $g\in\hol(U,\Omega)$ such that $g(z_j)=z_j$ for $j=1,2$.  By hypothesis there is a rational Kobayashi extremal function $h\in\hol(\D,\Omega)$ for $z_1, z_2$, so that equations \eqref{proph} hold.  Since $h$ is bounded it has no pole on $\T$, and therefore $h$ is analytic on some open neighbourhood of $\Delta$.  Hence there exists $t>1$ such that $h(t\D)\subset U$.  Let $h_t(\la)= h(t\la)$ for $\la\in\D$: then $g\circ h_t \in \hol(\D,\Omega)$ and $g\circ h_t(\la_j/t) = z_j$ for $j=1,2$.  Hence
\begin{align*}
\de_\Omega(z_1,z_2) &\leq \rho(\la_1/t, \la_2/t) \\
		&<\rho(\la_1,\la_2) \\
		&=\de_\Omega(z_1,z_2),
\end{align*}
which is a contradiction.  Thus $\id_\Omega$ is $2$-extremal. \hfill $\Box$

\subsection{Universal Carath\'eodory sets and $n$-extremal maps}\label{universal}
In this subsection we consider a question related to the ``Lempert duality'' between $\hol(\D,\Omega)$ and $\hol(\Omega,\D)$, mentioned in the Introduction, and to solvability criteria for interpolation problems.  Can we test a map $h\in\hol(\Omega_1,\Omega_2)$ for $n$-extremality by examining all its compositions with a suitable subset of $\hol(\Omega_2,\D)$?  By Proposition \ref{compose}, if $F\circ h$ is $n$-extremal for some $F\in\hol(\Omega_2,\D)$, then $h$ is $n$-extremal.  We ask whether there is a converse implication.

\begin{definition}\label{Caraset}
We say that a subset ${\rm C} \subset \hol(\Omega, \D)$ is a {\em universal Carath\'eodory set} for a domain $\Omega$ if, for every pair $z_1, z_2 \in \Omega$ there is a Carath\'eodory extremal function for $z_1, z_2$ that belongs to ${\rm C}$.
\end{definition}

\begin{example} {\rm For  many classical domains $\Omega$ there are small universal Carath\'eodory sets for $\Omega$. 

(i) If $\Omega = \D^d$ then the set of the $d$ co-ordinate functions is a  universal Carath\'eodory set for $\D^d$. 

(ii) For the ball $\mathbb{B}_d$ in $\C^d$, the projections onto the planes through the centre constitute a universal Carath\'eodory set. 

(iii) If $\Omega = \G$, the symmetrised bidisc, there is a $1$-parameter set $\{ \Phi_{\omega}: \omega \in \T \}$  (see Definition \ref{defPhi} below) that constitutes  a  universal Carath\'eodory set for $\G$ \cite{AY04}. 
}
\end{example}

Given a domain $\Omega$, a universal Carath\'eodory set $\mathrm{C}$ for $\Omega$ and an integer $n$, we may pose:
\begin{question}\label{Caraset-ext}  Is it true that $h\in\hol(\D,\Omega)$ is $n$-extremal if and only if $F \circ h$ is $n$-extremal for every $F \in {\mathrm C}$?
\end{question}

\begin{example}\label{exam-Cara-ext}{\rm The answer to Question \ref{Caraset-ext} depends on the domain $\Omega$ and on $n$.

(i) It is yes for the polydisc $\D^d$. It is easy to see that
$h =(h^1, \dots, h^d)$ is $n$-extremal if and only if some component $h^j$ is  $n$-extremal. 

(ii) If $\Omega = \G$, for $n \ge 3$, the answer to this question is no while for $n=2$ the answer is yes. In Proposition \ref{E_nu-munis-E_nu-1} we construct an analytic disc $h$ in $\G$ such that $h$ is $3$-extremal and yet, for all  $\omega \in \T$, $\Phi_{\omega} \circ h$ is not  $3$-extremal. 
}
\end{example}

\section{The symmetrised bidisc $\Gamma$} \label{symmetrised}
We began the study of the open symmetrised bidisc $\G$ in \cite{AY1} to \cite{AY06} with the aim of solving a special case of the $\mu$-synthesis problem of $H^\infty$ control: see our concluding reflections in Section \ref{conclud} below.  The original goal has still not been attained, but significant progress has been made, and the function theory of $\G$ has turned out to be of great interest to specialists in several complex variables.  The present study of $n$-extremal functions in $\hol(\D,\G)$ throws further light on interpolation problems for $\G$.

Here we summarise the relevant facts about $\G$.  We repeat the definitions:
 the {\em open} and  {\em closed symmetrised bidiscs} are 
 defined to be the sets
\begin{align*}
\G &\df\{(z+w,zw): |z|< 1, |w|< 1\} \quad \mbox{ and } \\
\Gamma &\df\{(z+w,zw): |z|\leq 1, |w|\leq 1\}
\end{align*}
respectively.  It is evident that the domain $\G$ is closely related to the bidisc, but $\G$ has a richer structure.  Its distinguished boundary is topologically a M\"obius band, and so is inhomogeneous (unlike that of the bidisc).
The (equal) invariant distances $\de_\G$ and $C_\G$ are less simple than for $\D^2$, but they can be calculated fairly explicitly \cite{AY04}.  The complex geodesics can also be described explicitly \cite{AY06,PZ}; they are rational of degree at most $2$, and so, by Proposition \ref{ident2ext}, $\id_\G$ is $2$-extremal.

The  distinguished boundary of $\G$ (or $\Gamma$) will be denoted by $b\Gamma$.  Thus $b\Gamma$
 is the \v{S}ilov boundary of the algebra of continuous functions on $\Gamma$ that are analytic in $\G$.  It is the symmetrisation of the 2-torus:
$$
b\Gamma= \{ (z+w,zw): |z|=|w|=1\}.
$$

Certain simple rational functions play a central role in the study of $\Gamma$.  
\begin{definition}\label{defPhi}
The function $\Phi$ is defined for $(z, s, p) \in \mathbb{C}^3$ such that $zs\neq 2$ by
$$
\Phi(z, s, p) = \frac{2zp -s}{2-zs}.
$$
We shall write $\Phi_z(s, p)$ as a synonym for $\Phi(z, s, p)$. 
\end{definition}
 In particular, $\Phi$ is defined and analytic on $\mathbb{D} \times
\Gamma$ (since $|s| \le 2$ when $(s, p) \in \Gamma$).   See \cite{AY2} for an account of how $\Phi$ arises from operator-theoretic considerations.

It will be useful to have criteria for a point of $\C^2$ to belong to $\Gamma$, $b \Gamma$  or the topological boundary $\partial\Gamma$.

\begin{proposition}\label{critGam} Let $(s,p) \in \C^2 $. Then

\noindent {\rm (1)} $(s,p) \in \G$ $\Leftrightarrow$ $|s- \bar{s} p| < 1 - |p|^2$; 

\noindent {\rm (2)} $(s,p) \in \Gamma$

\hspace{2cm} $\Leftrightarrow$ $|s|\le 2$ and  $|s- \bar{s} p|\le  1 - |p|^2$ 

\hspace{2cm} $\Leftrightarrow$ $|s| \le 2$ and, for all $\omega$ in a dense subset of $\T$, $|\Phi(\omega, s, p)| \le 1$;

\noindent {\rm (3)} $(s,p) \in b \Gamma$ $\Leftrightarrow$
$|s|\le 2$, $|p|=1$ and $s = \bar{s} p$;

\noindent {\rm (4)} $(s,p) \in \partial \Gamma$ 

\hspace{2cm} $\Leftrightarrow$ $|s|\le 2$ and  $|s- \bar{s} p| =  1 - |p|^2$ 

\hspace{2cm} $\Leftrightarrow$ there exist $z\in\T$ and $w\in\Delta$ such that $s=z+w, \ p=zw.$

Furthermore, for $\omega\in\T$ and $(s,p)\in\Gamma$,
\[
|\Phi_\omega(s,p)| = 1 \quad\mbox{ if and only if }\quad\omega(s-\bar s p)= 1-|p|^2.
\]
\end{proposition}
\begin{proof}
These statements are mainly contained in  \cite[Theorem 2.1 and Corollary 2.2]{AY04} and \cite[Introduction]{AY06}, {\em inter alia}, but statement (2) is a slight refinement.  The introduction of the dense subset of $\omega\in\T$ is occasioned by the fact that $\Phi_\omega$ is defined everywhere on $\Gamma$ except at the point $(2\bar\omega, \bar\omega^2)$.    

 To prove the second equivalence in (2), the fact that $(s,p)\in\Gamma$ implies that  $|s| \leq 2$ and $|\Phi_\omega(s,p)|\leq 1$ for all but at most one $\omega\in\T$ is in \cite[Theorem 2.1 and Corollary 2.2]{AY04}.  Conversely, if $|s|\leq 2$ and $|\Phi_\omega(s,p)|\leq 1$ for all  $\omega$ in a dense subset of $\T$, then, by continuity, for all $ \omega \in \T$, we have
\[
|2\omega p-s|^2 \leq |2-\omega s|^2.
\]
On expanding this inequality we find that
\[
4|p|^2 -4 \re ( \omega \bar{s} p) + |s|^2 \le 4 - 4 \re ( \omega s)+ |s|^2.
\]
Thus, for all $ \omega \in \T$,
\[
 \re ( \omega (s -\bar{s} p)) \le 1-|p|^2
\]
and therefore
\[
 |s -\bar{s} p| \le 1-|p|^2.
\]
It follows by the first equivalence in (2) that $(s,p) \in \Gamma$. 

The same calculation shows that, for fixed $\omega\in\T$ and $(s,p)\in\Gamma$,
\begin{align*}
|\Phi_\omega(s,p)|=1  &\Leftrightarrow \re (\omega(s-\bar s p)) = 1-|p|^2,
\end{align*}
and since, by part (2), $|s-\bar s p|\leq 1-|p|^2$, the last equation is true if and only if $\omega(s-\bar s p)= 1-|p|^2$.
\end{proof}

The variety
\[
\V \df \{(2z, z^2): z\in \C \}
\]
plays a special role in the study of $\Gamma$.  For one thing, 
$\V\cap \G$ is the orbit of $\{(0,0)\}$ under the automorphism group of $\G$.  We call $\V$ the {\em royal variety}.

\section{Interpolation in $\hol(\D,\Gamma)$ and a conjecture} \label{interp}

The (finite) interpolation problem for  $\hol(\D,\Gamma)$ is the following:

\noindent   {\em Given $n$ distinct points $\la_1, \dots, \la_n$ in the open unit disc $\D$ and $n$ points $z_1, \dots , z_n$ in $\Gamma$, find if possible an analytic function 
\beq\label{mainProb}
h:\D\to \Gamma \mbox{ such that } h(\la_j)=z_j \mbox{ for }j=1,\dots,n.
\eeq}

If $\Gamma$ is replaced by the closed unit disc $\Delta$ then we obtain the classical Nevanlinna-Pick problem \cite{Pick},  for which there is an extensive theory that furnishes among many other things a simple criterion for the existence of a solution $h$ and an elegant parametrisation of all solutions when they exist (see for example \cite{W,bgr,AgMcC}).  The classical results extend with appropriate modifications to a narrow class of other target sets, for example to the closed unit ball of the space of $k\times k$ matrices, or more generally, to closures of Cartan domains.  These extensions have applications in electrical engineering.  It would be useful for engineers if we could solve the finite interpolation problem for certain further sets, and a test case that has attracted much interest is the above problem of interpolation from $\D$ into $\Gamma$.

There is a satisfactory analytic theory of the problem \eqref{mainProb} in the case that the number of interpolation points $n$ is $2$ (a summary and references can be found in \cite{NJY11}), but we are still far from understanding the problem for a general $n\in\N$.  Here we introduce a sequence of necessary conditions  for the solvability of an  $n$-point $\Gamma$-interpolation problem and put forward a conjecture about sufficiency.  In Section \ref{inequations-for-E} we prove that these conditions are of strictly increasing strength.

  Consideration of some examples, including the case $n=2$,  led us to the following:

\begin{conjecture} {\bf The $\Gamma$-interpolation conjecture.}\label{conj3}
The $\Gamma$-interpolation data  
\[
\la_j \mapsto (s_j,\ p_j), \quad  1 \le j \le n,
\] 
  are  solvable  if and only if, for every Blaschke product $\up$ of degree
 at most $n-2$, the data
$$
\la_j \mapsto \frac{2\up(\la_j)p_j - s_j}{2-\up(\la_j)s_j},
\quad 1 \le j \le n,
$$
are solvable for the classical Nevanlinna-Pick problem.
\end{conjecture}
Here we say that
\begin{equation}\label{NPdata}
\lambda_j \mapsto z_j, \qquad 1 \le j  \le n,
\end{equation}
are {\em $\Gamma$-interpolation data} if $\lambda_1, \dots,
\lambda_n$ are distinct points in $\mathbb{D}$ and $z_1, \dots
z_n \in \Gamma$.  The data are {\em solvable} if there exists an analytic
function $f:\D\to \Gamma$  such that $f(\lambda_j) = z_j$ for
$j=1, \dots, n$; any such function is said to be a {\em solution} of
the $\Gamma$-interpolation problem \eqref{mainProb}  with data (\ref{NPdata}).    Observe that Pick's Theorem gives us an easily-checked criterion for the solvability of a Nevanlinna-Pick problem (see Proposition \ref{uppick}  below).

Conjecture \ref{conj3} is true in the case $n=2$  \cite{AY04}; see also \cite{NJY11}. 
We have no evidence for $n\ge3$ and we are open minded as to whether
or not it is likely to be true for all $n$.

We shall formalise the condition which appears in Conjecture \ref{conj3} and which plays an important role in the paper.
\begin{definition}\label{C_nu} 
Corresponding to interpolation data
\begin{equation}
\label{Gdata}
\la=(\la_1,\dots, \la_n), \quad z=(z_1, \dots, z_n),
\end{equation}
where $\lambda_1, \dots, \lambda_n$ are distinct points in $\D$
and $z_j = (s_j,p_j) \in \G$ for $j=1,\dots,n$,
we introduce:
\begin{center} \bf  Condition  $\mathcal{C}_\nu(\la,z)$ \end{center}
\noindent For every Blaschke product $\up$ of degree at most $\nu$, the Nevanlinna-Pick data
\begin{equation}
\label{upsdata}
\la_j \mapsto \frac{2\up(\la_j)p_j -s_j}{2-\up(\la_j)s_j}, \quad j=1,\dots,n,
\end{equation}
are solvable.
\end{definition}

Thus Conjecture \ref{conj3} can be stated: {\em Condition  $\mathcal{C}_{n-2}$ is necessary and sufficient for the solvability of an $n$-point $\Gamma$-interpolation problem.}

The conditions $\mathcal{C}_\nu$ are all necessary for the solvability of a $\Gamma$-interpolation problem.

\begin{theorem}\label{necG}
Let $\lambda_1, \dots, \lambda_n$ be distinct points in
$\mathbb{D}$ and let $z_j \in \G$ for $j=1, 2, \dots,
n$.  If there exists an analytic function $h:\mathbb{D} \to
\Gamma$ such that $h(\lambda_j) = z_j$ for $j=1, 2,
\dots, n$ then, for any function $\ups$ in the Schur class,
the Nevanlinna-Pick data \eqref{upsdata} are solvable. 
In particular, the condition
$\mathcal{C}_\nu(\la,z)$ 
holds for every non-negative integer $\nu$.
\end{theorem}

\begin{proof}  Suppose that the analytic function
$h$ exists as described.  Choose any function $\upsilon$ in the Schur class.  Then, by Proposition \ref{critGam}, for all $\lambda \in
\mathbb{D}$,
$$
|\Phi(\upsilon(\lambda), h(\lambda))|
\le 1.
$$
The function $g= \Phi \circ (\up, h)$ 
is analytic and bounded by 1 in $\mathbb{D}$, and 
satisfies the interpolation conditions (\ref{upsdata}).
Hence the Nevanlinna-Pick data (\ref{upsdata}) are indeed solvable.
In particular, the conclusion holds if $\upsilon$ is a Blaschke product,
and so condition $\mathcal{C}_\nu(\la,z)$ holds for every non-negative
integer $\nu$. \end{proof}

There is a special case  in which Condition $\mathcal{C}_0$ is sufficient as well as necessary.

\begin{theorem}
\label{flat}
Let $\lambda_1, \dots, \lambda_n$ be distinct points in $
\mathbb{D}$ and let $z_j=(s_j, p_j) \in \G, \  1 \le j \le n$.
If condition $\mathcal{C}_0(\la,z)$ holds and
the Nevanlinna-Pick
problem with data $\lambda_j \mapsto p_j$ is {\em extremally}
solvable then
$$
\la_j \mapsto z_j, \quad 1 \le j \le n,
$$
are solvable $\Gamma$-interpolation data.
\end{theorem}
This result is \cite[Theorem 5.2]{AY04}.

Condition $\mathcal{C}_0$ does not suffice for general
3-point interpolation problems, as will follow from Theorem \ref{C_nu-stronger-C_nu-1}. 

 Pick's Theorem enables us to recast the necessary condition
$\mathcal{C}_\nu$ in Theorem \ref{necG}
as the positivity of a pencil of matrices.
\begin{proposition}
\label{uppick}
If
$$
\la_j \mapsto z_j=(s_j,p_j), \quad 1 \le j \le n,
$$
are interpolation data for $\Gamma$ then condition 
$\mathcal{C}_\nu(\la_1,\dots,\la_n,z_1,\dots,z_n)$ holds
if and only if,
for every Blaschke product $\up$ of degree at most $\nu$,
\begin{equation}
\label{pick}
\left[\frac{1-\up(\lambda_i) p_i \bar p_j \overline
\up(\lambda_j) - \frac{1}{2} \up(\lambda_i) (s_i - p_i \bar s_j) -
\frac{1}{2} (\bar s_j - \bar p_j s_i) \overline \up
(\lambda_j)
   -\tfrac 14   (1-\ups(\la_i) \bar\ups(\la_j) ) s_i \bar s_j
}{1-\lambda_i \overline \lambda_j}\right]^n_{i, j=1}
\ge 0.
\end{equation}
\end{proposition}

\begin{proof}
By Pick's Theorem \cite[Theorem 1.3]{AgMcC}, for any function $\up$ in the
Schur class, the Nevanlinna-Pick data
(\ref{upsdata}) are solvable if and only if
\begin{equation}
\label{hpick}
\left[\frac{1- \bar w_i w_j)}{1-\bar\lambda_i \lambda_j}\right]^n_{i, j=1} \ge 0,
\end{equation}
where
\begin{eqnarray*}
w_j
&=& \frac{2\up(\lambda_j) p_j -s_j}{2-\up(\lambda_j)s_j}.
\end{eqnarray*}
On conjugating the inequality (\ref{hpick}) by
$$
\mathrm{diag} \{2 - \up(\lambda_1) s_1, \dots,
2-\up(\lambda_n)s_n\}
$$
we deduce that the inequality (\ref{hpick}) holds if and only if
$$
\left[\frac{(2-\up(\lambda_i)s_i)(2-\overline
\up(\lambda_j)\bar s_j) - (2\up(\lambda_i) p_i -s_i)(2\overline
\up(\lambda_j)\bar p_j - \bar s_j)}{1-\lambda_i\overline
\lambda_j} \right]^n_{i, j=1} \ge 0,
$$
which simplifies to the given condition (\ref{pick}).
Hence condition $\mathcal{C}_\nu(\la,z)$ holds if and only if,
for every Blaschke product $\up$ of degree at most $\nu$,
the inequality (\ref{pick}) holds.
\end{proof}

An analogue of condition $\mathcal{C}_0$ for interpolation into the symmetrised polydisc in $\C^N$, for $N \ge 2$, was found by D. Ogle \cite[Corollary 5.2.2]{ogle}. However, when $N \ge 3$, this necessary condition is insufficient even for two-point interpolation problems \cite[Observation 1.3]{bharali}.

\section{Extremality in condition $\mathcal{C}_\nu$} \label{extremality}

To prove that condition $\mathcal{C}_\nu$ suffices for the solvability of an $n$-point Nevanlinna-Pick problem for $\Gamma$ it is enough to prove it in the case that $\mathcal{C}_\nu$ holds {\em extremally}.  Let us make this notion precise.

Recall that $\Gamma$-interpolation data $\la_j \mapsto z_j, \ 1\leq j\leq n,$ are defined to satisfy condition $\mathcal{C}_\nu$ if, for every 
Blaschke product $\up \in \B l_\nu$ of degree at most $\nu$, the data
\begin{equation}
\label{upsolve}
\la_j \mapsto \Phi(\up(\la_j), z_j), \quad 1 \le j \le n,
\end{equation}
are solvable for the classical Nevanlinna-Pick problem.
If, in addition, there exists $m \in  \B l_\nu$ such that the data 
\[
\la_j \mapsto \Phi(m(\la_j), z_j), \quad 1 \le j \le n,
\]
are {\em extremally} solvable Nevanlinna-Pick data, then we shall say that the data $\la_j \mapsto z_j, \ 1\leq j\leq n,$ {\em satisfy $\mathcal{C}_\nu$ extremally}, or the condition $\mathcal{C}_\nu(\la, z)$ {\em holds extremally}. Here $\la=(\la_1,\dots, \la_n)$ and $ z=(z_1, \dots, z_n)$.

It is well known (e.g. \cite{AgMcC}) that  Pick's criterion for the solvability of a classical Nevanlinna-Pick
problem is expressible by an operator norm inequality; hence condition $\mathcal{C}_\nu$ can be expressed this way.   Let
\begin{equation}\label{defM}
\mathcal{M} = \mathrm{span}\ \{K_{\lambda_1}, \dots,
K_{\lambda_n}\} 
\subset H^2,
\end{equation}
where $K$ is the Szeg\H{o} kernel. Consider $\Gamma$-interpolation data $\la_j \mapsto z_j, \ 1\leq j\leq n$, and 
introduce, for any function $\up$ in the Schur class, the operator
$X(\up)$ on $\mathcal{M}$ given by
\beq\label{defX}
X(\up)K_{\la_j}= \overline{ \Phi(\up(\la_j),z_j)}K_{\la_j}, \qquad 1\leq j\leq n.
\eeq
Pick's Theorem, as reformulated by Sarason \cite{sar}, asserts that the Nevanlinna-Pick data
\begin{equation}
\label{Pick-data}
\la_j \mapsto \Phi(\up(\la_j),z_j), \qquad 1\leq j\leq n,
\eeq
are solvable if and only if the operator
$X(\up)$ is a contraction. Furthermore, the Nevanlinna-Pick data \eqref{Pick-data} are extremally solvable if and only if $\|X(\up) \|=1$.

 Thus $\mathcal{C}_\nu(\la, z)$ holds if and only if
\begin{equation}
\label{Xule1}
  \sup_{\up \in \B l_{\nu}} \|X(\up) \| \leq 1.
\end{equation}
\begin{proposition} \label{equivCnu}
For any $\Gamma$-interpolation data $\la_j \mapsto z_j, \ 1\leq j\leq n$, and $\nu \ge 0$, the following conditions are equivalent.

{\rm (i)} $\mathcal{C}_\nu(\la, z)$  holds extremally;

{\rm (ii)} $\sup_{\up \in \B l_{\nu}} \|X(\upsilon)\| = 1;$

{\rm (iii)} $\mathcal{C}_\nu(\la, z)$  holds and there exist
$ m \in \B l_\nu$ and $q \in \B l_{n -1}$ such that 
\begin{equation}
\label{Phiq}
\Phi(m(\la_j), z_j)= q(\la_j), \quad j =1, \dots, n,
\end{equation}

Moreover, when condition {\rm (iii)} is satisfied for some $m \in \B l_\nu$, there is a {\em unique}  $q \in \B l_{n -1}$ such that equations \eqref{Phiq} hold. If, furthermore, the 
$\Gamma$-interpolation data $\la_j \mapsto z_j, \ 1\leq j\leq n$, are solvable by an analytic function $h=(s,p):\D\to\Gamma$, then
\beq\label{propq}
\frac{2mp - s}{2 - ms} = q.
\eeq
\end{proposition}
\begin{proof} (i) $\Rightarrow$ (ii). Assume (i): then the inequality \eqref{Xule1} holds. Furthermore there exists $m \in \B l_\nu$ such that 
\begin{equation}
\label{Phi-data1}
\la_j \mapsto \Phi(m(\la_j), z_j), \quad 1 \le j \le n,
\end{equation}
are extremally solvable  Nevanlinna-Pick data, which implies that 
$$
\sup_{\up \in \B l_{\nu}} \|X(\upsilon)\| = 1.
$$
Hence (ii) holds.

 (ii) $\Rightarrow$ (iii). Suppose (ii). Take a maximizing sequence $(\up_k) \in \B l_\nu$ for this supremum: it has a
locally uniformly convergent subsequence, whose limit is a Blaschke product $m$ of degree $d(m)$ at most $\nu$. 
We have $\|X(m)\| =1$, so that the Nevanlinna-Pick data 
\begin{equation}
\label{Phi-data2}
\la_j \mapsto \Phi(m(\la_j), z_j), \quad 1 \le j \le n,
\end{equation}
are extremally solvable.  Hence there is a unique interpolating
function, $q$ say, in the Schur class,  
and $q \in \B l_{n -1}$  \cite{W}. Thus (iii) holds.

 (iii) $\Rightarrow$ (i). Let (iii) hold. Since 
the Nevanlinna-Pick data 
\begin{equation}
\label{Phim}
\la_j \mapsto \Phi(m(\la_j), z_j), \quad 1 \leq j \leq n,
\end{equation}
are  solvable by a Blaschke product of degree $\le n-1$, they are extremally solvable \cite[Lemma 6.19 and Theorem 1.3]{AgMcC}. Thus (i) holds.

Hence (i), (ii) and (iii) are equivalent.

If   $m \in \B l_\nu$ is such that \eqref{Phiq} holds for some  $q \in \B l_{n -1}$  then  the 
Nevanlinna-Pick problem with data \eqref{Phim}, being extremally solvable, has a unique solution in the Schur class. The solution is necessarily $q$
\cite[Theorem 6.4]{AgMcC}.

Finally, if $h=(s,p):\D\to\Gamma$ is a solution of the $\Gamma$-interpolation problem
$\la_j \mapsto z_j, \ 1\leq j\leq n$, then 
$\Phi \circ (m, h) \in \Schur$
solves the problem \eqref{Phim}, and so, by uniqueness, equals $q$.
That is $$
\frac{2mp - s}{2 - ms} = q.  
$$
\end{proof}

Suppose that $\mathcal{C}_\nu$ holds for $n$-point interpolation data $\la_j\mapsto (s_j,p_j), \ 1\leq j\leq n$.  On replacing $\la_j$ by $r\la_j$ for a suitable $r\in (0,1]$ we can arrange that $\mathcal{C}_\nu$ holds extremally for the modified data.  If condition $\mathcal{C}_\nu$ suffices for solvability of an $n$-point interpolation problem for $\Gamma$ in the extremal case then there exists an analytic function $g:\D\to\Gamma$ such that $g(r\la_j)=z_j,\ 1\leq j\leq n$.  Then the function $h(\la)=g(r\la)$ solves the original problem.  This justifies the claim at the beginning of this section that it suffices to prove sufficiency in Conjecture \ref{conj3} for the case  that $\mathcal{C}_\nu$ holds extremally.

We shall say that any Blaschke product $m$ with the properties described in Proposition \ref{equivCnu}(iii)  is an {\em auxiliary extremal} for the condition $\mathcal{C}_\nu(\la,z)$.

Let us consider the degrees of auxiliary extremals $m$ associated with data $\la_j\mapsto (s_j,p_j), \ 1\leq j\leq 3,$ that satisfy $\mathcal{C}_1$ extremally.  It is far from the case that $m$ is uniquely determined, or even that the degree $d(m)$ is unique for a particular set of data. 

\begin{examples} \label{exdm}    \rm
Let $\la_1, \, \la_2, \, \la_3$ be any three distinct points in $\D$ and let $ 0 < r < 1$.   In each of the following examples $h$ is an analytic function from $\D$ to $\G$ and the data $\la_j \mapsto h(\la_j), \, 1\leq j\leq 3,$ satisfy $\mathcal{C}_1$ extremally. \\
(1)  Let $h(\la)=(2r\la,\la^2)$.  Every degree 0 inner function $m \in\T$ is an auxiliary extremal for $\mathcal{C}_1$; there is no auxiliary extremal of degree $1$. \\
(2)  Let $h(\la)= (r(1+\la), \la)$.  Every $m \in \B l_1$  is an auxiliary extremal for $\mathcal{C}_1$.  The corresponding $q$ has degree  $d(m)+1$. \\
(3)  Let 
\[
h(\la)= \left(2(1-r)\frac{\la^2}{1+r\la^3}, \frac{\la(\la^3+r)}{1+r\la^3}\right), \qquad \la\in\D.
\]
  The function $m(\la)=-\la$ is an auxiliary extremal for $\mathcal{C}_1$; there is no auxiliary extremal of degree $0$.  Here $q(\la)=-\la^2$.
See Proposition \ref{E_nu-munis-E_nu-1} for a more general example and Theorem \ref{C_nu-stronger-C_nu-1} for applications.\\
(4) Let $f$ be a Blaschke product of degree $1$ or $2$ and let $h=(2f,f^2)$.  Every  $m \in \B l_1$ is an auxiliary extremal and, for every $m$, we have $q=-f$.
\end{examples}

\section{$\Gamma$-inner functions}\label{Gamin}


\begin{definition}\label{Gam-in-funct} A  {\em $\Gamma$-inner function} is an analytic function $h : \D \to \Gamma$ such that the radial limit
\begin{equation}\label{radial}
\lim_{r \to 1-} h(r \lambda) \in b \Gamma
\end{equation}
for almost all $\lambda \in \T$.
\end{definition}

\rm
By Fatou's Theorem, the radial limit (\ref{radial}) exists for almost all 
 $\lambda \in \T$ with respect to Lebesgue measure. 
Observe that, in view of Proposition \ref{critGam}(3), if $h= (h_1, h_2)$ is a
 $\Gamma$-inner function, then $h_2$ is an inner function on $\D$ in the conventional sense.

The easiest way to construct a $\Gamma$-inner function is to symmetrise a pair of inner functions.

\begin{example}\label{boring} {\rm Let $\phi$ and $\psi$ be inner functions on $\D$. Then
\begin{equation}\label{phipsi}
h = (\phi + \psi, \phi \psi)
\end{equation}
is $\Gamma$-inner. In particular, $h = (2 \phi, \phi^2)$ is 
$\Gamma$-inner; this example has the property that $h(\D)$ lies in the royal variety $\V$.
}
\end{example}

The analysis of the  $\Gamma$-inner function \eqref{phipsi} reduces to the study of pairs of conventional inner functions. 

\begin{proposition}\label{critBoring} Let $h =(s,p)$ be a $\Gamma$-inner function. Then $h$ is of the form $(\phi + \psi, \phi \psi)$ for some pair 
 $\phi$, $\psi$ of  inner functions on $\D$ if and only if $s^2 -4p$ has an analytic square root on  $\D$. 
\end{proposition}  

\begin{proof} Necessity is trivial, since $s^2 -4p= (\phi - \psi)^2$.
Conversely, if $s^2 -4p= g^2$ where $g$ is analytic on $\D$, then we may take $\phi= \frac{1}{2}(s -g)$ and $ \psi= \frac{1}{2}(s +g)$; we find that $ \phi, \  \psi$ are inner and $h = (\phi + \psi, \phi \psi)$.
\end{proof}

\begin{example}\label{flatGeo} Let $|\beta|< 1$. The function 
\begin{equation}\label{betaGeo}
h(\la) = (\beta \la + \bar{\beta}, \la)
\end{equation}
is $\Gamma$-inner and is not the symmetrisation of a pair of inner functions.
\end{example}
\begin{proof} Write $h =(s,p)$. It is easy to see from Proposition \ref{critGam} that $h$ is $\Gamma$-inner. We claim that $s^2 -4p$ does not have an analytic square root on $\D$; to prove this it suffices to show that $s^2 -4p$ has a {\em simple} zero in $\D$. Now
\[
s^2 -4p = (\beta \la + \bar{\beta})^2 -4  \la= \beta^2 \la^2 +2 (|\beta|^2 -2)\la +(\bar{\beta})^2.
\]
If $\beta =0$ then  $s^2 -4p$ has a simple zero at $0$.
Otherwise, $s^2 -4p$ is a quadratic with discriminant
\[
4[(|\beta|^2 -2)^2 - |\beta|^4] =16 (1 - |\beta|^2 )
\]
which is nonzero for $|\beta| <1$, and so $s^2 -4p$ has two distinct zeros. The sum and product of these zeros satisfy
\[
|{\rm sum}|= 2 \frac{2- |\beta|^2}{|\beta|^2} >2, \;\; 
|{\rm product}|= \left|\frac{\bar{\beta}^2}{\beta^2} \right|=1.
\]
The former inequality implies that the zeros of $s^2 -4p$ do not both lie on $\T$ and the latter then shows that one lies in $\D$ and one in $\C \setminus \Delta$. Consequently  $s^2 -4p$ has a simple zero in $\D$, and so $h$ is not the  symmetrisation of a pair of inner functions.
\end{proof}

The function $h$ given by equation \eqref{betaGeo} is not only $\Gamma$-inner -- it is a {\em complex geodesic} of $\G$. In the case of the function $h$ from  equation \eqref{betaGeo} the simplest left inverse is of course the projection $(s,p) \mapsto p$. The domain $\G$ also has complex geodesics of degree 2 \cite[Theorem 0.2]{AY06}: these afford a further class of interesting  $\Gamma$-inner functions.

\subsection{New $\Gamma$-inner functions from old} \label{newOld}

There are three easy ways to construct new $\Gamma$-inner functions from a given $\Gamma$-inner function.

\begin{proposition} \label{3ways}
Let $h =(s,p)$ be a $\Gamma$-inner function.

{\rm (i)} For any inner function $\phi$, $h \circ \phi$ is $\Gamma$-inner.

{\rm (ii)} For any inner function $\up$, $(\up s, \up^2 p)$ is $\Gamma$-inner.

{\rm (iii)} If
\begin{equation}\label{rmin}
0 \le r \le \min \left\{\frac{2}{\|s\|_\infty}, \inf_{\D} \frac{1- |p|^2}{|s - \bar{s}p|} \right\}
\end{equation}
then $(rs,p)$ is $\Gamma$-inner. In particular, $(rs,p)$ is $\Gamma$-inner
 whenever $0 \le r \le 1$.
\end{proposition}

 All three statements follow easily from Proposition \ref{critGam}. The minimum on the right hand side of inequality \eqref{rmin} can be greater than $1$: for $h(\la) =(\beta \la +\bar{\beta}, \la)$ in Example \ref{flatGeo}, the minimum is $\frac{1}{|\beta|}$. 

Also under the heading of ``new from old" comes the intriguing fact that the $\Gamma$-inner functions have a non-obvious semigroup structure; see \cite[Corollary 16(1a)]{Ma} and  \cite[Theorem 3.4]{Cost05}.\\

\begin{proposition} \label{semigp} If $(s,p)$ and $(t,q)$ are $\Gamma$-inner functions then so is $(\frac{1}{2}st,pq)$.
\end{proposition}
\begin{proof} Clearly $|\tfrac{1}{2}st | \le 2$ on $\D$, $pq$ is inner and
\[
\tfrac{1}{2}st= \tfrac{1}{2}(\bar{s}p)(\bar{t}q)=\left(\tfrac{1}{2}\bar{s}\bar{t}\right) pq
\]
on $\T$.
\end{proof}

Note that Proposition \ref{3ways}(ii) is a special case of Proposition \ref{semigp}, since $(2 \up, \up^2)$ is obviously  $\Gamma$-inner for
 any inner function $\up$.
The constant function $(2,1)$ is an identity for the semigroup of $\Gamma$-inner functions. The only idempotents of the semigroup are $(2,1)$ and $(0,1)$. The only elements that have inverses in the semigroup are the constant functions of the form $(2\omega,\omega^2)$ for some $\omega \in \T$.
The semigroup structure will not play a role in this paper.

\subsection{Rational $\Gamma$-inner functions} \label{RatGamIn}

The following example shows that {\em not every pole of $p$ need be a pole of $s$} for a rational $\Gamma$-inner function $(s,p)$.

\begin{example} \label{surprise}{\rm Let $a \in \D \setminus  \{ 0\}$ and let
\[
h(\la) = (s(\la),p(\la)) = \left(\frac{c \la}{1-\bar{a} \la}, \frac{\la(\la-a)}{1-\bar{a} \la} \right)
\]
for some $c \in \R$ such that $|c| \le 2(1- |a|)$.
It is easy to check that $h$ is $\Gamma$-inner, and whereas $p$ has a pole at infinity, $s(\infty)= -\frac{c}{\bar{a}} \neq \infty$. Similarly, if $b \in  \D \setminus \{0 \}$ and $(s^{\flat}, p^{\flat}) = h \circ B_b$, then 
$(s^{\flat}, p^{\flat})$ is  $\Gamma$-inner, $ p^{\flat}$ has a pole at 
$B_{-b}(\infty)= \frac{1}{\bar{b}}$ and $s^{\flat}(\frac{1}{\bar{b}})= -\frac{c}{\bar{a}}\neq \infty$.
}\qed
\end{example}

In Example \ref{surprise}, the pole $\frac{1}{\bar{b}}$ of $ p^{\flat}$ that is not a pole of  $ s^{\flat}$ corresponds to the common zero $b$ of $ s^{\flat}$ and $ p^{\flat}$. This phenomenon is general.

\begin{proposition} \label{polesZeros} Let $(s,p)$ be a rational $\Gamma$-inner function.

{\rm (i)} If $a \in \C \cup \{\infty\}$ is a pole of $s$, of multiplicity $k \ge 1$, then $a$ is a pole of $p$ of multiplicity at least $k$.

{\rm (ii)} If $a \in \C \cup \{\infty\}$ is a pole of $p$, of multiplicity $k \ge 0$, and $\frac{1}{\bar{a}}$ is a zero of $s$ of multiplicity $\ell \ge 0$, then $a$ is a pole of $s$ of multiplicity at least $k- \ell$.
\end{proposition}

\begin{proof} (i) The equations
\[
s(\la) = \overline{s(\la)}p(\la)= s^\vee(\bar{\la}) p(\la)=s^\vee \left(\tfrac{1}{\la} \right) p(\la)
\]
hold for all $\la \in \T$.
Since the first and last terms are rational functions,
\begin{equation} \label{scup}
s(\la)=s^\vee \left(\tfrac{1}{\la} \right) p(\la) \;\;\;\text{for all } \; \la \in \C.
\end{equation}
Hence, if $a \in \C$,
\begin{equation} \label{poleq}
(\la -a)^{k -1}s(\la)=s^{\vee}\left(\tfrac{1}{\la} \right)(\la -a)^{k -1} p(\la) \;\;\;\text{for all } \; \la \in \C.
\end{equation}
Since $a$ is a pole of $s$, $|a| >1$, hence $\left|\frac{1}{a}\right| <1$,
and so 
$s^{\vee}$ is analytic at $\frac{1}{a}$. On letting $\la \to a$ in equation \ref{poleq} we find that $(\la -a)^{k -1} p(\la) \to \infty$ as $\la \to a$. Thus $a$ is a pole of $p$ of  multiplicity at least $k$.

Now suppose that $\infty$ is a pole of $s$ of multiplicity  $k \ge 1$.
Then $0$ is  a pole of $s\left(\tfrac{1}{\la} \right)$ of multiplicity  $k$,
so that $\la^{k -1} s\left(\tfrac{1}{\la} \right) \to \infty$ as $\la \to 0$. From the relation \eqref{scup} we have 
\[
\la^{k -1} s\left(\tfrac{1}{\la} \right) = s^\vee(\la) \la^{k -1} p\left(\tfrac{1}{\la} \right).
\]
Since $s^\vee$ is analytic at $0$, we have
$\lim_{\la \to 0} s^\vee(\la) = s^\vee(0)$. Therefore
\[
\lim_{\la \to 0}  \la^{k -1} p\left(\tfrac{1}{\la} \right) = \infty,
\]
which is to say that $p$ has a pole of  multiplicity at least $k$ at $\infty$.

(ii) Again the statement follows on letting $\la \to a$ in equation \eqref{scup}, since  $s^\vee \left(\tfrac{1}{\la} \right)$ has a zero of  multiplicity  $\ell$ at $a$.

Now suppose that $\infty$ is a pole of $p$ of multiplicity  $k \ge 1$ and $0$ is a zero of $s$ of  multiplicity  $\ell \ge 0$; then $\infty$ is a pole of $s$ of  multiplicity $k- \ell$.
From the relation \eqref{scup}, for all $\la \in \C \setminus \{0\}$, we have 
\[
s \left(\tfrac{1}{\la} \right)= s^\vee(\la) p\left(\tfrac{1}{\la} \right),
\]
and so
\begin{equation} \label{s-p-poles}
\frac{\la^{k -1}}{\la^\ell} s\left(\tfrac{1}{\la} \right) = \frac{s^\vee(\la) }{\la^\ell} \la^{k -1} p\left(\tfrac{1}{\la} \right).
\end{equation}
Since $s^\vee$ is analytic at $0$ and has  a zero at $0$  of  multiplicity  $\ell \ge 0$, we have
\[
\lim_{\la \to 0} \frac{s^\vee(\la) }{\la^\ell} = c \neq 0 \; \text{where}\; c \in \C.
\]
By assumption, $p(\la)$ has a pole of  multiplicity  $k$ at $\infty$, and so 
$
\lim_{\la \to 0}  \la^{k -1} p\left(\tfrac{1}{\la} \right) = \infty
$.
Therefore, by equation \ref{s-p-poles},
\[
\lim_{\la \to 0}  \la^{k -\ell-1} s\left(\tfrac{1}{\la} \right) = \infty.
\]
which is to say that  $s\left(\tfrac{1}{\la} \right)$ has a pole of  multiplicity  $k-\ell$ at $0$. Thus $s(\la)$ has a pole of  multiplicity  $k-\ell$ at $\infty$.
\end{proof}

\begin{remark} {\em In Proposition \ref{polesZeros} 
we allow the possibility that $\ell > k$, in which case $a$ is a zero of $s$ of  multiplicity $\ell -k$. In Example \ref{flatGeo} we consider the  rational $\Gamma$-inner function
$h(\la) =(\beta \la +\bar{\beta}, \la)$. The function  $s(\la) =\beta \la +\bar{\beta}$ has a zero of  multiplicity $\ell =1$ at $ \la = - \frac{\bar{\beta}}{\beta}$. The function  $p(\la) =\la$ has a pole of  multiplicity $k =0$ at $ \la = - \frac{\bar{\beta}}{\beta}$. Thus we have $\ell > k$.
}
\end{remark}
\begin{corollary} \label{denomsp} If $(s,p)$ is a rational $\Gamma$-inner function
then $s$ and $p$ can be written as ratios of polynomials with the same denominators. More precisely, let
\begin{equation} \label{formp}
p(\lambda) = c \frac{\lambda^k \tilde{D}_p(\lambda)}{D_p(\lambda)}
\end{equation}
where $|c|=1$, $k \ge 0$, $D_p$ is a polynomial of degree $n$ such that $D_p(0)=1$ and  $\tilde{D}_p(\lambda) =\lambda^n D_p^\vee \left(\tfrac{1}{\la} \right) $. Then $s$ is expressible in the form 
\begin{equation} \label{forms}
s(\lambda) = \frac{\lambda^{\ell} N_s(\lambda)}{D_p(\lambda)} 
\end{equation}
where $0 \le \ell \le \frac{1}{2}(n +k)= \frac{1}{2} d(p)$,  and $N_s$ is a polynomial of degree $ d(p) - 2 \ell$ such that $N_s(0) \neq 0$. Moreover, if 
$N_s(\la) = \sum_{j =0}^{n+k -2\ell} b_j \la^j$  then 
\begin{equation} \label{relbj}
b_j = c \bar{b}_{n+k-2 \ell -j} \;\; \text{for } \; j =0,1, \dots, n+k-2 \ell.
\end{equation}
The degree of $s$ is at most $\max \{n+k-\ell , n \}.$
\end{corollary}
\begin{proof}
Every finite Blaschke product is expressible in the form \eqref{formp}.
Since every pole of $s$ is a pole of $p$, $s$ can be written as a  ratio of polynomials with  denominator $D_p$, though not necessarily in its lowest terms (cf. Example \ref{surprise}).

Let $\ell \ge 0$ be the multiplicity with which $s$ vanishes at $0$: then $s$ can be written in the form \eqref{forms} for some polynomial $N_s$.
From the relation $s(\la)=s^\vee \left(\tfrac{1}{\la} \right) p(\la)$ for all $\la \in \C$ we deduce
\begin{equation} \label{form-s-p}
\frac{\lambda^{\ell} N_s(\lambda)}{D_p(\lambda)}=\frac{ 1}{\lambda^{\ell}}
\frac{ N^\vee_s(1/\lambda)}{D^\vee_p(1/\lambda)}
 c \lambda^k \frac{\lambda^n D^\vee_p(1/\lambda)}{D_p(\lambda)}.
\end{equation}
Hence,  for all $\la \in \C$,
\[
 N_s(\lambda) = c \lambda^{n +k -2 \ell} N^\vee_s(1/\lambda).
\]
That is, if $N_s$ has degree $d$, 
\[
N_s(\la) = \sum_{j =0}^{d} b_j \la^j= c \lambda^{n +k -2 \ell}\sum_{i =0}^{d} \frac{\bar{b_i}}{ \la^i}.
\]
The constant term $b_0$ is nonzero and is the term of lowest degree, and hence
$ 0 = n +k -2 \ell -d$ and
\[
b_0 = c \bar{b}_d = c  \bar{b}_{n +k -2 \ell}.
\]
Similarly, on equating coefficients of $\la^j$, we obtain equation \eqref{relbj}. From the fact that $d \ge 0$ we conclude that $2 \ell  \le n +k$.
\end{proof}

\section{The classes $\E_{\nu k}$}\label{Eclasses}

Proposition \ref{equivCnu} tells us that if $h \in {\rm Hol}(\D, \Gamma)$ and $\la_1, \dots,\la_n$ are distinct points in $\D$, then the 
$\Gamma$-interpolation data $\lambda_j \mapsto h(\la_j)$ satisfy $C_{\nu}(\lambda, h(\la))$ {\em extremally} if and only if there exists $m \in \B l_{\nu}$ such that 
$\Phi \circ (m, h) \in \B l_{n-1}$. This leads us to introduce the following classes of rational  $\Gamma$-inner functions.

\begin{definition}\label{E_nu_k}
For $\nu\geq 0, \ k\geq 1$ we say that the function $h $ is in  $\E_{\nu k}$ 
if $h=(s,p)\in\hol(\D,\Gamma)$ is rational and there exists $m \in \B l_{\nu}$ such that 
\[
\frac{2 m p -s}{2- m s}  \in \B l_{k-1}.
\]
\end{definition}

\begin{remark} {\rm
It is obvious that, for every $\nu \ge 0$,
\[
\E_{\nu 1} \subset \E_{\nu 2} \subset \dots  \subset \E_{\nu k}\subset \E_{\nu, k+1}\subset \dots, 
\]
and, for every $k \ge 1$,
\[
\E_{0 k} \subset \E_{1 k} \subset \dots  \subset \E_{\nu k}\subset \E_{\nu+1, k}\subset \dots. 
\]
We shall see in Section \ref{kExtremEnuk} that there is a strong connection between the class $\E_{\nu k}$ and $k$-extremality.  In fact if $h\in\E_{\nu k}$ and $h(\D)$ meets $\G$ then $h$ is $k$-extremal, while if the $\Gamma$-interpolation Conjecture is true, then every $k$-extremal in $\hol(\D,\Gamma)$ belongs to $\E_{k-2,k}$ (Theorem \ref{Enukextrem} and Observation \ref{C-n-extrem}).
}
\end{remark}
It is not obvious that the functions in $\E_{\nu n}$ are $\Gamma$-inner, but it is so.
\begin{theorem}\label{ratGin}
Let $h\in\hol(\D,\Gamma)$.  If there exists an inner function $m$ such that $\Phi\circ (m,h)$ is inner then $h$ is $\Gamma$-inner.  If furthermore $h\in\E_{\nu n}$ for some integers $\nu\geq 0, n\geq 1$ then $h$ is rational of degree at most $2n-2$.
\end{theorem}
\begin{proof}
Let $h=(s,p)$ and let  $q=\Phi\circ (m,h)$.   Consider first the case that $mq=1$ identically.  Then $m, q$ are constant -- say $m=\omega, \ q=\bar\omega$ for some $\omega\in\T$.  The relation $\Phi(\omega,s,p)=\bar\omega$ tells us that $p$ is constant and equal to $\bar\omega^2$, and then the fact that
\[
|s-\bar s p| \leq 1-|p|^2 =0
\]
shows that $\omega s = \bar\omega \bar s$, and hence $\omega s$ is a real constant.  Thus in this case $h$ is constant and equal to $(x\bar\omega, \bar\omega^2)$ for some $x\in [-2,2]$, which implies that $h$ is $\Gamma$-inner and rational of degree $0$.

The remaining case is that $mq$ is not identically equal to $1$, and so $1-mq \neq 0$ a.e. on $\T$.
Since $q$ is inner, for almost all $\la\in\T$ we have 
\[
|\Phi_{m(\la)}(s(\la),p(\la))|=1, 
\]
and so, by Proposition  \ref{critGam},
\beq\label{msp}
m(\la)(s-\bar s p)(\la) = 1-|p(\la)|^2.
\eeq
On rearranging the equation $\Phi(m,s,p)=q$ we obtain
\beq\label{propsmpq}
s=2\frac{mp-q}{1-mq},
\eeq
and so equation \eqref{msp} becomes
\[
m\left(2\frac{mp-q}{1-mq}-2\frac{\bar m\bar p-\bar q}{1-\bar m\bar q} p\right)(\la) = 1-|p(\la)|^2 \quad \mbox{ a.e. on } \T.
\]
Since $|mq(\la)| = 1$ a.e. we may multiply above and below in the second term on the left hand side to obtain
\begin{align*}
1-|p(\la)|^2 &= 2m\left(\frac{mp-q}{1-mq}-\frac{q\bar p - m}{mq-1}p\right)(\la) \\
	&= \frac{2m}{1-mq}(mp-q-(q\bar p -m)p)(\la) \\
	&= -\frac{2m}{1-mq}(\la)(1-|p(\la)|^2) \qquad \mbox{ a.e.}
\end{align*}
Thus 
\[
(1-|p(\la)|^2) \frac{1+mq}{1-mq}(\la) = 0
\]
a.e. on $\T$, and so $p$ is inner.

Since $h(\D)\subset \Gamma$, we have $|s(\la)| \leq 2$ a.e. on $\T$, and by equation \eqref{msp} we also have
\[
s(\la) = \overline{ s(\la)}p(\la) \quad \mbox{ a.e. }
\]
and so, by Proposition \ref{critGam}(3), $h(\la)\in b\Gamma$ for almost all $\la\in\T$, that is, $h$ is $\Gamma$-inner.

We prove the second statement.  Suppose $h\in\E_{\nu n}$.  There exist $m\in\B l_\nu, \ q\in\B l_{n-1}$ such that $\Phi\circ (m,h) =q$.  From the relation \eqref{propsmpq} and the fact that $|s(\la)|\leq 2$, we have
\[
|mp-q|^2 \leq |1-mq|^2
\]
on $\T$. Since $m, p$ and $q$ are all unimodular on $\T$, on expanding we find that
\[
1 -2 \re (mp \bar{q}) +1 \le 1 - 2 \re (m q) +1,
\]
on $\T$, and so
\[
\re (mq) \leq \re (mp \bar{q})
\]
everywhere on $\T$.  Since $mp(\la), \ mp\bar q(\la) \in\T$ for all $\la\in\T$, it follows that if $mp\bar q(\la) =-1$ then also $mq(\la) = -1$.  Now the finite Blaschke product $mq$ takes on the value $-1$ precisely $d(mq)$ times on $\T$, and hence $mp\bar q$ can take the value $-1$ at most $d(mq)$ times.  The winding number wno of $mp\bar q$ about $0$ is thus at most $d(mq)$, that is
\[
 d(m)+d(p)-d(q) = \mathrm{wno} (mp\bar q) \leq d(m)+d(q),
\]
and therefore
\[
d(p) \leq 2d(q) \leq 2n-2. 
\]
\end{proof}

\begin{proposition} If
$h \in \E_{\nu n}$ then, for any choice of distinct points $\lambda_1, \dots, \lambda_n \in \D$, the data
\[
\lambda_j \mapsto h(\lambda_j), \; 1 \le j \le n,
\]
satisfy $\mathcal{C}_{\nu}(\lambda, z)$ extremally where $z_j = h(\lambda_j)$, $ 1 \le j \le n$.
\end{proposition}
\begin{proof} It follows from Proposition \ref{equivCnu}.
\end{proof}

\subsection{Phasar derivatives}\label{phasar}
To study the classes $\mathcal{E}_{\nu k}$ we shall need the following basic notions relating to functions on the unit circle.

\begin{definition}\label{phasar-derivative}
 For any differentiable function $f: \T \to \C \setminus \{0\}$
 the {\em phasar derivative} of $f$ at $z= \e^{\ii \ta} \in \T$
is  the derivative with respect to $\ta$ of the argument of $f(\e^{\ii \ta})$ at $\ta$; we denote it by $Af(z)$.
\end{definition}
Thus if $f(\e^{\ii \ta}) = R(\ta) \e^{\ii g(\ta)}$ is differentiable   where $g(\ta) \in\R$ and $R(\ta) >0$ then $g$ is differentiable on $[0,2\pi)$ and the phasar derivative of $f$ at $z= \e^{\ii \ta} \in \T$ is equal to
\begin{equation}\label{phasar-deriv-arg}
Af(\e^{\ii \ta}) = \frac{d}{d \ta} \arg f(\e^{\ii \ta}) = g'(\ta).
\end{equation}

Clearly, for differentiable functions $\psi, \phi : \T \to \C \setminus \{0\}$ and
for any $c \in \C \setminus \{0\}$, we have
\beq\label {phasar-deriv-sum}
A(\psi \phi ) =A \psi +A \phi \; \quad \text{and} \; \quad A(c \psi  ) =A \psi.
\eeq

The following is a simple calculation.
\begin{proposition}\label{AofInner}  Let $\phi: \T \to \C \setminus \{0\}$ be a rational inner function. Then, for all $\lambda \in \T$,
\begin{equation} \label{phasar-deriv-inner}
 A\phi(\lambda)= \lambda \frac{\phi'(\lambda)}{\phi(\lambda)}.
\end{equation}
\end{proposition} 

\begin{proposition} \label{phasar-p-Blaschke} 
{\rm (i)} Let 
$$
B_\alpha(z) = \frac{z-\alpha}{1-\overline \alpha z}
$$
be a Blaschke factor for  $\alpha \in \D$. Then the phasar derivative
$A B_\alpha >0$ on  $\T$.

{\rm (ii)} Let $p$ be a rational inner function. Then
the phasar derivative $A p (\lambda) >0$ for  all $\lambda \in \T$.
\end{proposition}

\begin{proof} (i) By Proposition \ref{AofInner},  for all $\alpha \in \D$ and $\lambda \in \T$, we have
\beq\label{phasar-Bl}
 A B_\alpha (\lambda) = 
\lambda \frac{B'_\alpha(\lambda)}{B_\alpha(\lambda)}= 
\frac{1-|\alpha|^2}{|1-\overline{\alpha} \lambda|^2} >0.
\eeq

(ii) For such $p$, it is well known that there exist $c \in \T$ and $
\alpha_1, \alpha_2, \dots,\alpha_n \in \D$ such that
$p = c B_{\alpha_1}  \dots B_{\alpha_n}$.
Then, by Remark \ref{phasar-deriv-sum}  and Part (i),
\[
A p  = A B_{\alpha_1}  + A B_{\alpha_2}  +  \dots + A B_{\alpha_n}  >0.
\]
\end{proof}

\subsection{Cancellations and the classes $\E_{\nu k}$}\label{cancels}
Cancellations in the functions $\Phi\circ (\up, h)$, where $h$ is a rational $\Gamma$-inner function, are at the heart of the technical results of this paper.

Consider a rational $\Gamma$-inner function $h=(s,p) $ where
\[
s(\lambda) = \frac{\lambda^{\ell} N_s(\lambda)}{D_p(\lambda)} \;\;\text{and}\; \;
p(\lambda) = \frac{\lambda^k \tilde{D}_p(\lambda)}{D_p(\lambda)}
\]
where $\ell \le \frac{1}{2} d(p)$ and $d(N_s) = d(p) - 2 \ell$; see
Corollary \ref{denomsp}.

For  $\nu \ge 0$ and $\up \in \B l_{\nu}$, consider the function
\[
\Phi \circ (\up, h) = \frac{2 \up p -s}{2- \up s}. 
\]
What is the degree of the rational function $\Phi \circ (\up, h)$?
Since the denominator of $s$ divides the denominator of $p$ (see Corollary \ref{denomsp}), the function $\Phi \circ (\up, h)$ is a Blaschke product of degree at most $ d(\up p)$. If $\Phi \circ (\up, h)$ has no cancellations then it has degree exactly $d(\up p)$.   It transpires that cancellations can only happen at special points on the unit circle.

\begin{definition} A point $\la \in \Delta$ is a {\em royal node} of  a $\Gamma$-inner function $h$ if and only if  $h(\la)$ is in the royal variety $\V = \{(2z, z^2): z \in \C \}$.
\end{definition}

Clearly a point  $\la \in \Delta$ is a royal node of  a $\Gamma$-inner function $h=(s,p)$ if and only if 
\beq \label{royal-node} 
s(\la)^2= 4p(\la).
\eeq

\begin{definition} A $\Gamma$-inner function $h =(s,p)$ is {\em full} if $\|s \|_{\infty} =2$.
\end{definition}

\begin{lemma}\label{royal=full} Let $h$ be a rational $\Gamma$-inner function. 

{\rm (i)} The royal nodes of $h =(s,p)$ on $\T$ are precisely the points $\la \in \T$ such that $|s(\la) | =2$.

{\rm (ii)}  The function $h$ has a royal node  on $\T$ if and only if $h$ is full.
\end{lemma}
\begin{proof} (i) By Remark \ref{royal-node},
$\la \in \T$ is a royal node of  $h=(s,p)$ if and only if $s^2(\la)= 4p(\la)$. Note that
\[
|s(\la) | =2 \quad \Leftrightarrow  \quad  s(\la)\overline{s(\la)} =4.
\]
By Proposition \ref{critGam}(3), for every $\la \in \T$,
 $\overline{s(\la)}p(\la)=s(\la)$. Hence $|s(\la) | =2$ implies $s(\la)\overline{s(\la)}p(\la) =4p(\la)$, and so $s^2(\la)= 4p(\la)$.
On the other hand, by Proposition \ref{critGam}(3), for every $\la \in \T$,
$|p(\la)|=1$. Therefore, $s^2(\la)= 4p(\la)$ implies $|s(\la)|^2 = 4 |p(\la)|=4$.

(ii) If $h$ has a royal node $\la$ on $\T$ then, by Part (i) $h$ is full.
Suppose $\|s \|_{\infty} =2$. By the Maximum Principle, there is a $\la \in \T$ such that $|s(\la) | =2$. By Part (i), $\la \in \T$ is a royal node of  $h$.
\end{proof}

\begin{proposition} \label{canc=royal} Let $h =(s,p)$ be a rational $\Gamma$-inner function.  Suppose $\up$ is a finite Blaschke product such that 
\[
\Phi \circ (\up, h) = \frac{2 \up p -s}{2- \up s} 
\]
has a cancellation at a point $\zeta\in\C$. Then $h$ is full, $\zeta \in \T$, $\zeta$ is a royal node for $h$, $\up(\zeta)=\frac{1}{2}\overline{ s(\zeta)}$  and $|s(\zeta)|=2$. 
\end{proposition}
\begin{proof}  By assumption, $\Phi \circ (\up, h)$ has a cancellation at $\zeta$, and so
\[
(2 \up p -s)(\zeta) = 0 = (2- \up s)(\zeta).
\]
Thus $\up(\zeta) s(\zeta)=2$ and $2 \up(\zeta)s(\zeta) p(\zeta) =s^2(\zeta)$. Therefore, $s^2(\zeta)= 4p(\zeta)$, and so $\zeta$ is a royal node for $h$. One can also see that $2 \up^2(\zeta) p(\zeta) =\up(\zeta) s(\zeta)=2$. Since $\up$ is a finite Blaschke product, the equality $\up^2(\zeta) p(\zeta) =1$ implies that $\zeta \in \T$.  By Lemma \ref{royal=full},  $h$ is full and $|s(\zeta)|=2$. 
Note that
$\up(\zeta) s(\zeta)=2$ and $|s(\zeta)|=2$ imply that $\up(\zeta)=\frac{1}{2}\overline{ s(\zeta)}$.
\end{proof}

\begin{theorem} \label{crit-canc=royal} Let $h =(s,p)$ be a nonconstant rational $\Gamma$-inner function and let  $\up$ be a finite Blaschke product.
Then   
$\Phi \circ (\up, h) $
has a cancellation at $\zeta$ if and only if the following conditions are satisfied: $\zeta \in \T$, $\zeta$ is a royal node for $h$ and  $\up(\zeta)=\frac{1}{2}\overline{ s(\zeta)}$.
Moreover $\Phi \circ (\up, h)$ has at most one cancellation at any royal node $\zeta$.
\end{theorem}
\begin{proof} Necessity follows by Proposition \ref{canc=royal}.

Conversely, let $\zeta$ be a royal node for $h$ on $\T$, let
$h(\zeta) = ( 2 \bar{\omega}, \bar{\omega}^2)$ where $\omega \in \T$ and let $\up(\zeta) = \omega$. Hence
$
(2 \up p -s)(\zeta) = 2 \omega \bar{\omega}^2 -2 \bar{\omega}=0
$
and 
$
(2- \up s)(\zeta) = 2 - \omega \cdot 2\bar{\omega}= 0.
$
Thus $\Phi \circ (\up, h)$ has at least one  cancellation at $\zeta$. 

Suppose there are 2 cancellations at $\zeta$, so that 
\[
(2 \up p -s)'(\zeta) = 0 = (2- \up s)'(\zeta).
\]
Then
\begin{eqnarray}
\label {deriv-1}
(2 \up p -s)'(\zeta) 
&=& 2 \up'(\zeta) p(\zeta) +2 \up(\zeta) p'(\zeta) -s'(\zeta)\nonumber\\
&=& 2 \up'(\zeta) \bar{\omega}^2  +2 \omega p'(\zeta) -s'(\zeta)
 = 0
\end{eqnarray}
and
\begin{eqnarray}
\label {deriv-2}
(2- \up s)'(\zeta) 
&=& -  \up'(\zeta) s(\zeta) - \up(\zeta) s'(\zeta)\nonumber\\
&=& - \up'(\zeta) 2 \bar{\omega} -  \omega s'(\zeta)
 = 0.
\end{eqnarray}
Hence we have $\up'(\zeta)  = -\frac{1}{2}  \omega^2 s'(\zeta)$.
Then equation \eqref{deriv-1} can be written as
\[
 - 2 \frac{1}{2}  \omega^2 s'(\zeta)\bar{\omega}^2  +2 \omega p'(\zeta) -s'(\zeta) = 0.
\]
Therefore $ s'(\zeta) = \omega p'(\zeta)$.
By Proposition \ref{AofInner}, the phasar derivatives
\[
Ap(\zeta)= \zeta \frac{p'(\zeta)}{p(\zeta)}= \zeta  \frac{\bar{\omega}s'(\zeta)}{\bar{\omega}^2} = \zeta \omega s'(\zeta).
\]
and
\[
A\up(\zeta)= \zeta \frac{\up'(\zeta)}{\up(\zeta)}= \zeta  \frac{ -\frac{1}{2}  \omega^2 s'(\zeta) }{\omega} = -\frac{1}{2} \zeta \omega s'(\zeta) =-\frac{1}{2}Ap(\zeta).
\]
In view of Proposition \ref{phasar-p-Blaschke}, we have a contradiction since the phasar derivatives $Ap(\zeta) >0$ and $A\up(\zeta) \ge 0$.
Thus the function $\Phi \circ (\up, h)$ has exactly one  cancellation at $\zeta$.
\end{proof}

\begin{proposition}\label{h-not-full-E} Let $h =(s,p)$ be a rational  $\Gamma$-inner function. Suppose $h$ is not full. Then, for all $\nu \ge 0$,
$h \in \E_{\nu, \nu+ d(p) +1} \setminus \E_{\nu, \nu + d(p)}$.
\end{proposition}
\begin{proof}
Since  $h$ is not full,  $|s(\la) | < 2 $ for all $\la \in \Delta$. By Proposition \ref{canc=royal}, for all $\nu \ge 0$ and all $\up \in \B l_{\nu}$, $\Phi \circ (\up, h)$ 
has no cancellations, and so has degree $d(\up p)$.
\end{proof}

\begin{corollary}\label{h=(0,p)-E} Let $h$ be a rational  $\Gamma$-inner function of the form $(0,p)$.  Then, for all $\nu \ge 0$,
$h \in \E_{\nu, \nu+ d(p) +1} \setminus \E_{\nu, \nu + d(p)}$.
\end{corollary}
\begin{proof} It follows from Proposition \ref{h-not-full-E}.
\end{proof}

\begin{corollary}\label{Npoints-canc-royal=1}
Let $h =(s,p)$ be a  rational  $\Gamma$-inner function, and let
$\zeta_1, \dots, \zeta_N$ be distinct  royal nodes for $h$ on $\T$. Let
$h(\zeta_j) = ( 2 \bar{\omega_j}, \bar{\omega_j}^2)$ for $j =1, \dots, N$. 
If $\up$ is a finite Blaschke product such that $\up(\zeta_j)= \omega_j$  for $j =1, \dots, N$, then 
$\Phi \circ (\up, h)$ is a Blaschke product of degree $d(\up p) -N$.
\end{corollary}
\begin{proof} By Theorem \ref{crit-canc=royal}, the function $\Phi \circ (\up, h)$ has one cancellation at each point $\zeta_j$, $j =1, \dots, N$.
\end{proof}

For a given rational $\Gamma$-inner $h=(s,p)$, what can we expect of the $\E$-classes to which $h$ belongs, in terms of the degree $d(p)$? For every $\nu \ge 0$, we trivially have $h \in \E_{\nu, \nu+ d(p) +1}$; the interesting question is whether $h \in \E_{\nu k}$ for some $k$ less than 
$\nu+ d(p) +1$. By Proposition \ref{h-not-full-E}, if $h$ is not full the answer is no. If $h$ is  full we can always arrange one cancellation at a royal node, and so $h \in \E_{\nu, \nu+ d(p)}$ for every $\nu \ge 0$.
If $h$ has $N$ royal nodes in $\T$ we can arrange $N$ cancellations in $\Phi \circ (m, h)$ by choosing $\nu$ sufficiently large that there exists $m \in \B l_{\nu}$ that maps the royal nodes of $h$ in $\T$ to the required target points in $\T$, as in Corollary \ref{Npoints-canc-royal=1}, and then we shall have 
$h \in \E_{\nu, \nu+ d(p)- N +1}$. However, the question as to how large 
$\nu$ must be is subtle, as we shall see in the following two examples; they show that, even for the symmetrization of a pair of Blaschke products, it is a delicate issue whether one can achieve $3$ cancellations in $\Phi \circ (\up, h)$ with $\up \in \B l_1$.

\begin{example}\label{E-1-k+2} Let $\psi$ be a rational inner function on $\D$.  The rational $\Gamma$-inner function
\begin{equation}\label{h_psi}
h_\psi(\la) =\left( \la + \la \psi(\la),\la^{2} \psi(\la) \right), \;\la \in \D,
\end{equation}
lies in $ \E_{1, d(\psi)+2} \setminus \E_{1, d(\psi)+1}$.
\end{example}
\begin{proof} 
By Theorem \ref{crit-canc=royal}, for any Blaschke product $\up$,  the function $\Phi \circ (\up, h_\psi)$
has a cancellation at $\zeta$ if and only if  $\zeta \in \T$, $\zeta$ is a royal node for $h_\psi$ and  $\up(\zeta)=\frac{1}{2}\overline{ s_\psi(\zeta)}$. The royal nodes for $h_\psi= (s_\psi, p_\psi)$ are the roots of
$s_\psi^2(\la)-4p_\psi(\la)=0$, that is, $(\la - \la \psi(\la))^2=0$.
Therefore the royal nodes $\omega_j$ for $h_\psi$ on $\T$ are the  roots of $\psi(\la)=1$, that is, $\omega_j \in \T$ such that $\psi(\omega_j)=1$,
 $\; j= 1, \dots, d(\psi)$.
For $j= 1, \dots,d(\psi)$, we have $\tfrac{1}{2}\overline{s_\psi(\omega_j)}= \tfrac{1}{2} \overline{(\omega_j + \omega_j\psi(\omega_j))} = \overline{\omega_j}$.

To prove that  $h_\psi$ is not in  $ \E_{1, d(\psi)+1}$ we must show that, for all $m \in \B l_1$, the rational function $\Phi \circ (m,h_\psi)$ is not in $\B l_{d(\psi)}$. It is enough to show that, for all $m \in \B l_1$, the rational function $\Phi \circ (m,h_\psi)$ cannot have $3$ cancellations.
To get  $3$ cancellations we need $3$ royal nodes $\omega_{j_1}$,
$\omega_{j_2}$ and $\omega_{j_3}$, $1 \le j_1 < j_2 <j_3 \le d(\psi)$, with
$m (\omega_{j_1}) =\overline{\omega_{j_1}}$, $m (\omega_{j_2}) =\overline{\omega_{j_2}}$ and $m (\omega_{j_3}) =\overline{\omega_{j_3}}$. This is impossible since the points $\overline{\omega_{j_1}}$, $\overline{\omega_{j_2}}$, $\overline{\omega_{j_3}}$ are in the opposite cyclic order to 
$\omega_{j_1}$, $\omega_{j_2}$, $\omega_{j_3}$ on $\T$.

Let us show that $ h_\psi \in \E_{0, d(\psi)+2}$, that is, there is $\omega \in \T$ such that $\Phi \circ (\omega,h_\psi)  \in \B l_{d(\psi)+1}$. Take  $\omega= \tfrac{1}{2}\overline{s_\psi(\omega_1)}$;  then, by Theorem \ref{crit-canc=royal},
 $\Phi \circ (\omega,h_\psi)$ has exactly one cancellation at $\omega_1$. Therefore $\Phi \circ (\omega,h_\psi)$
has degree $d(\psi)+1$, and so $\Phi \circ (\omega,h_\psi) \in \B l_{d(\psi)+1}$. 
Note that $\E_{0, d(\psi)+2} \subset \E_{1, d(\psi)+2}$.
Therefore  $ h_\psi \in \E_{1, d(\psi)+2} \setminus \E_{1, d(\psi)+1}$.
\end{proof}

The next example looks similar to $h_\psi$, and yet here one {\em can} achieve $3$ cancellations with $\up$ of degree $1$.
\begin{example}\label{E-1-2j+4}  For any positive integer $j$ the rational  $\Gamma$-inner function
\begin{equation}\label{h_2j+5}
h_j(\la) =\left( \la^2 + \la^{2j+3},\la^{2j+5} \right), \;\la \in \D,
\end{equation}
 belongs to  $\E_{1, 2j+4} \setminus \E_{0,  2j+4}$.
\end{example}
\begin{proof} By Example \ref{boring}, $h_j$ is $\Gamma$-inner.
By Theorem \ref{crit-canc=royal}, for any Blaschke product $\up$,  the function
\[
\Phi \circ (\up, h_j) = \frac{2 \up p_j -s_j}{2- \up s_j} 
\]
has a cancellation at $\zeta$ if and only if $\zeta \in \T$, $\zeta$ is a royal node for $h_j$ and $\up(\zeta)=\frac{1}{2}\overline{ s_j(\zeta)}$. 
The royal nodes for $h_j=(s_j, p_j)$ are the roots of
$s_j^2(\la)-4p_j(\la)=0$, that is, $(\la^2 - \la^{2j+3})^2=0$.
Therefore the royal nodes $\omega_k$ for $h_j$ on $\T$ are the $(2j +1)$th roots of $1$, that is, 
\[
\omega_k = e^{2 \ii \pi k/(2j+1)}, \; k= 0, \dots, 2j.
\]

To prove that  $h_j$ is not in  $ \E_{0,  2j+4}$ we must show that, for all $\omega  \in \T$, the rational function $\Phi \circ (\omega,h_j)$ is not in $\B l_{2j+3}$. If $\omega= \tfrac{1}{2}\overline{s_j(\omega_k)}$ for some $k= 0, \dots, 2j$,  then,
by Theorem \ref{crit-canc=royal}, $\Phi \circ (\omega,h_j)$ has exactly one cancellation at $\omega_k$. Therefore $\Phi \circ (\omega,h_j)$
has degree $2j +4$, and so $\Phi \circ (\omega,h_j)$ is not in $\B l_{2j+3}$. If $\omega \neq \tfrac{1}{2}\overline{s_j(\omega_k)}$ for all $k= 0, \dots, 2j$, then $\Phi \circ (\omega,h_j)$ has degree $2j +5$.
Thus $h_j$ is not in  $ \E_{0,  2j+4}$.

Let us show that $h_j \in \E_{1,  2j+4}$, that is, there is $m \in \B l_1$ such that 
\[
\Phi \circ (m,h_j)=\frac{2 m p_j -s_j}{2- m s_j}  \in \B l_{2j +3}.
\]
For all $k= 0, \dots, 2j$,  $\tfrac{1}{2}\overline{s_j(\omega_k)}= \tfrac{1}{2} \overline{(\omega_k^2 + \omega_k^{2j+3})} = \overline{\omega_k}^2$.
Note that the points $\omega_0=1,\; \omega_1, \; \omega_1^{j +1}$ have the same cyclic order as $\overline{\omega_0}^2=1,\; \overline{\omega_1}^2, \; \overline{\omega_1^{j +1}}^2= \overline{\omega_1}$ on $\T$.
Take $m \in \B l_1$ such that 
\[
m(\zeta) = \tfrac{1}{2}\overline{s_j(\zeta)}=\overline{\zeta}^2
\]
for $\zeta = \omega_0, \omega_1, \omega_1^{j +1}$.
By Corollary \ref{Npoints-canc-royal=1}, $\Phi \circ (m, h)$ has  cancellations at $\omega_0,\; \omega_1, \; \omega_1^{j +1}$. Thus  $\Phi \circ (m, h_j)$ is a Blaschke product of degree $d(p)+d(m) -3=  2j +3$.
Therefore $h_j \in \E_{1, 2j+4} \setminus \E_{0,  2j+4}$.
\end{proof}

\section{Superficial $\Gamma$-inner functions  and the classes $\E_{\nu 1}$}\label{superficial}

In the next three sections of the paper we derive some further information about the two-dimensional array of classes $(\mathcal{E}_{\nu k})$.

For any inner function $\ph$ and $\omega\in\T$ the function $h=(\omega+\ph,\omega\ph)$ is $\Gamma$-inner, and has the property that $h(\la)$ lies in the topological boundary $\partial\Gamma$ of $\Gamma$ for all $\la\in\D$.  We shall prove not only a converse of this statement, but also the fact that all the classes in the first column of the array $(\mathcal{E}_{\nu k})$ consist of precisely this type of $\Gamma$-inner function.
\begin{definition}\label{superfic} A function $h\in \hol(\D,\Gamma)$ is {\em superficial} if $h(\D) \subset \partial\Gamma$.
\end{definition}

\begin{lemma}\label{Phiconstant} Let  $h$ be a nonconstant $\Gamma$-inner function such that 
\[
\Phi_{\omega}\circ h =  \kappa
\]
identically on $\D$ for some $\omega, \kappa \in \T$. Then 
$\kappa = - \bar{\omega}$ and  $h = ( \omega p + \bar{\omega}, p)$ for some inner function $p$.
\end{lemma}
\begin{proof} Let $h=(s,p)$: then $p$ is a nonconstant inner function.  By assumption,
\[
\Phi_{\omega}(s, p) = \frac{2 \omega p -s}{2- \omega s} = \kappa.
\]
Hence
\[
s = 2 \frac{\omega p -\kappa}{1-\omega \kappa}.
\]
Since $|s| \le 2$, we have $|1-\omega \kappa| \ge |\omega p - \kappa|$.
Therefore for all $\lambda \in \Delta$, 
\[
|1-\omega \kappa| \ge |\omega p(\lambda) - \kappa|= | p(\lambda) - \bar{\omega} \kappa|.
\]
Since $p$ is nonconstant there is $\lambda_0 \in \T$ such that
$p(\lambda_0) = - \bar{\omega} \kappa$, and so $
|1-\omega \kappa| \ge 2$. Thus we have $\omega \kappa = -1$ and $\kappa = -\bar{\omega} $. Hence
\[
s = 2 \frac{\omega p -\kappa}{1-\omega \kappa}  = \omega p +\bar{\omega}.
\]
\end{proof}

\begin{proposition}\label{superficial-nonconst}  A $\Gamma$-inner function $h $ is superficial  if and only if there is an $\omega \in \T$ and an inner function $p$ such that $h=(\omega p + \bar{\omega},p) $.
\end{proposition}
\begin{proof} $\Leftarrow$   Since the topological boundary $\partial\Gamma$ of $\Gamma$ comprises the points $\{(z+w,zw): |z|=1, \ |w|\leq 1\}$, it follows that $h(\D) \subset \partial\Gamma$.

$\Rightarrow$ Consider a superficial $\Gamma$-inner $h=(s,p)$; necessarily $p$ is inner.  If $h$ is constant we may write $(s,p)=(\omega +z, \omega z)$ for some $\omega\in\T$ and $z\in\Delta$; then we have $z=\bar\omega p$ and $h=(\omega p + \bar{\omega},p) $.

Now suppose that $h$ is nonconstant.  It follows that $s$ is nonconstant, and hence $|s(0)| <2$. Observe that, for any $(s_0,p_0)\in\partial\Gamma$ we have $|s_0 - \bar s_0 p_0|=1-|p_0|^2$.  Now
$h(0) \in \partial \Gamma$, and so there exists $\omega \in \T$ such that
\[
\omega \left(s(0) -\overline{s(0)} p(0)\right)= |s(0)- \overline{s(0)}p(0)| = 1 -|p(0)|^2.
\]
A simple calculation (or \cite[Theorem 2.5]{AY04}) shows that
\[
\left|\Phi_{\omega}(s(0), p(0))\right|  = \left|\frac{2 \omega p(0) -s(0)}{2- \omega s(0)}\right| = 1.
\]
Therefore $\Phi_{\omega} \circ h$ is an inner function which takes a value of modulus $1$ at $0$, and so it is a constant function, with value $\kappa$, say. 
By Lemma \ref{Phiconstant}, $\kappa = - \bar{\omega}$ and $s= \omega p + \bar{\omega}$, and so $h = ( \omega p + \bar{\omega}, p)$.
\end{proof}

The image of a  function in $\hol(\D,\Gamma)$ is either contained in or disjoint from $\partial\Gamma$.
\begin{lemma} \label{notsuper}
If $h\in \hol(\D,\Gamma)$ is not superficial then $h(\D) \subset \G$.
\end{lemma}
\begin{proof}
Let $h=(s,p)$.  Suppose that $h(\D)$ is not contained in $\G$: then there exists $\la_0\in\D$ such that $h(\la_0) \in \partial\Gamma$.  Let $u(\la)$ be the spectral radius of the matrix
\[
\begin{pmatrix}  0 & 1 \\ -p(\la) & s(\la) \end{pmatrix}
\]
for $\la\in\D$.  By Vesentini's Theorem \cite[Theorem 2.3.32]{dales}, $u$ is subharmonic in $\D$.  Since $u(\la)$ is the maximum of the moduli of the roots of the equation $z^2-s(\la)z+p(\la) =0$, we have
\[
u(\la)= \max\{|z|,|w|\}  \mbox{ where } s(\la)=z+w, \, p(\la)= zw.
\]
Since $(s(\la),p(\la))\in\Gamma$, we have $0 \leq u(\la) \leq 1$ for all $\la \in \D$, and $u(\la) = 1$ if and only if $(s(\la),p(\la)) \in \partial \Gamma$.  By hypothesis $u$ attains its maximum at a point $\la_0\in\D$, and hence $u$ is constant and equal to $1$ on $\D$.  Consequently $h(\la) \in\partial \Gamma$ for all $\la\in\D$. 
\end{proof}

\begin{proposition}\label{E_0_1=superficial}  The class $\E_{0 1}$ consists of the superficial rational $\Gamma$-inner functions.
\end{proposition}
\begin{proof} By definition, a rational function $f=(s,p) \in \E_{0 1}$ if and only if there are $\omega, \kappa \in \T$ such that  
\begin{equation}\label{super-s}
\Phi_{\omega}\circ (s, p) =\frac{2 \omega p -s}{2- \omega s} = \kappa.
\end{equation}

First consider the case that $p$ is nonconstant. For
$f=(s,p) $, the conditions of Lemma \ref{Phiconstant} are satisfied. Hence,
we have $\kappa = - \bar{\omega}$, $s= \omega p + \bar{\omega}$ and  $f = ( \omega p + \bar{\omega}, p)$.
By Proposition \ref{superficial-nonconst}, $f$ is a superficial  $\Gamma$-inner function.

In the case that $p$ is constant, since $f=(s,p)$ is a $\Gamma$-inner function, by Proposition \ref{critGam}, $|s|\le 2$, $|p|=1$, $s = \bar{s} p$ and  $f=(s,p)\in b \Gamma$, and so is superficial.
\end{proof}

\begin{theorem}\label{E_nu_1=superficial} For every $\nu \ge 1$, the class $\E_{\nu 1}$ is equal to $\E_{0 1}$ and  consists of the superficial rational $\Gamma$-inner functions.
\end{theorem}
\begin{proof}
By Definition \ref{E_nu_k} and Theorem \ref{ratGin}, the function $h = (s,p) \in \E_{\nu 1}$ 
if $h=(s,p) $ is rational $\Gamma$-inner and there exists $m \in \B l_{\nu}$ such that 
\[
\frac{2 m p -s}{2- m s}= \kappa \in \T.
\]
Then
\[
s = 2 \frac{m p -\kappa}{1- m \kappa}.
\]
Since $|s| \le 2$ on $\T$, we have
\[
|m p - \kappa|^2 \le |1- m \kappa|^2.
\]
Since $m, p$ and $\kappa$ are all unimodular on $\T$, on expanding we find that
\[
1 -2 \re (mp \bar{\kappa}) +1 \le 1 - 2 \re (m \kappa) +1,
\]
on $\T$, and so
\begin{equation}\label{Re(mk)<=Re(mpk)}
\re (m \kappa)(\lambda) \le \re (mp \bar{\kappa})(\lambda)
\end{equation}
for all  $\lambda \in \T$.

Therefore, if  $m \kappa =1$ at some point $\lambda \in \T$ then 
$mp \bar{\kappa}=1$ at the point $\lambda$. Suppose the Blaschke product $m$ has degree $\nu' \le \nu$, so that  $m \kappa =1$ at $\nu'$ 
distinct points $\lambda_1, \dots, \lambda_{\nu'} \in \T$. Then also
$mp \bar{\kappa}=1$ at the points $\lambda_1, \dots, \lambda_{\nu'} \in \T$. Hence 
\[
(mp \bar{\kappa})(\lambda_j)= (m \kappa) (\lambda_j) (\bar{\kappa}^2 p )(\lambda_j)=\bar{\kappa}^2 p(\lambda_j)=1
\]
for $j =1, \dots, \nu'$. Thus $p=\kappa^2$ at 
 $\nu'$ distinct points $\lambda_1, \dots, \lambda_{\nu'} \in \T$, and so the degree $d(p) \ge d(m)= \nu'$.

We claim  that $d(m)=0$.
Suppose $d(m)>0$; then $d(p) \ge d(m)> 0$, and $p$ is nonconstant.
Since $p$ is nonconstant and rational, by Proposition \ref{phasar-p-Blaschke}, 
the phasar derivative $A (p)$ of $p$ is strictly positive on $\T$.
Observe that phasar derivatives satisfies the following inequality
\[
A(mp \bar{\kappa}) =A(mp) = A(m) +A(p) > A(m) = A(m \kappa)
\] 
on $\T$.  Recall that $|mp \bar{\kappa}|=1 =|m \kappa| =1$ on  $\T$, and, as we have shown above, there is a point $\lambda_0 \in \T$ such that
 $m \kappa(\lambda_0) =1$ and $mp \bar{\kappa}(\lambda_0) =1$.
 Therefore, at some $\lambda' \in \T$ close to the point $\lambda_0$, we have
\[
 \re (mp \bar{\kappa})(\lambda') < \re (m \kappa)(\lambda')
\]
which is a contradiction to the inequality \eqref{Re(mk)<=Re(mpk)}.
Hence $d(m)=0$, and therefore, $h = (s,p) \in \E_{0 1}$.
\end{proof}

\section{The classes $\E_{\nu k}$ and $k$-extremals}\label{kExtremEnuk}

In this section we show that the elements of $\E_{\nu k}$ are $k$-extremal for $\Gamma$ and that, if the $\Gamma$-interpolation Conjecture \ref{conj3} is true, then every $k$-extremal rational $\Gamma$-inner function belongs to $\E_{k-2,k}$.

For $\zeta \in\C$ and $(s,p)\in\C^2$ we define
\[
\zeta \cdot(s,p) = (\zeta s, \zeta^2 p) \quad \mbox{ and } \quad \zeta \cdot \Gamma = \{\zeta \cdot z: z\in\Gamma\}.
\]

\begin{theorem}\label{Enukextrem}
If $h \in \mathcal{E}_{\nu k}$, where $\nu\geq 0$ and $k\geq 2$, and $h$ is not superficial then $h$ is $k$-extremal for $\hol(\D,\Gamma)$.
\end{theorem}
\begin{proof}
Let $h\in\mathcal{E}_{\nu k}$ be not superficial.  There exist $m \in \B l_\nu$ and $q\in \B l_{k-1}$ such that $\Phi\circ(m,h) = q$.  Suppose that $h$ is not $k$-extremal: then there exist $k$ distinct points $\la_1,\dots,\la_k \in\D$, an $r_0 > 1$ and a function $f\in \hol(r_0\D, \Gamma)$ such that $f(\la_j)=h(\la_j)$ for $j=1,\dots,k$.  Since $h$ is not superficial, the points $h(\la_j) \in\G$, by Lemma \ref{notsuper}, and so, by the same lemma, $f(r_0\D) \subset \G$.

Pick any $r_1$ in the interval $(1,r_0)$: then $f(r_1\Delta)$ is a compact subset of $\G$.  Now 
\[
f(r_1\Delta) \subset \bigcup_{0<\rho<1} \,  \rho \cdot\G =\G,
\]
and hence there exists $\rho \in (0,1)$ such that $f(r_1\Delta)\subset \rho\cdot\G\subset\rho\cdot\Gamma$.

Observe that, for $\la\in\Delta$ and $(s,p) \in\Gamma$ we have
\[
\Phi(\la,\rho\cdot(s,p)) =\Phi(\la,\rho s,\rho^2 p) = \frac{2\la\rho^2p-\rho s}{2-\la\rho s} = \rho\Phi(\rho\la, s, p) \in \rho\Delta.
\]
Thus
\[
\Phi(\Delta \times \rho\cdot\Gamma) \subset \rho\Delta \subset \D.
\]
Furthermore, $\Phi$ is analytic on $(\rho^{-1}\D) \times \rho\cdot\Gamma$.  Hence, by continuity of $\Phi$ and compactness of $\rho\cdot\Gamma$, there is a neighbourhood $U$ of $\Delta$ such that
\[
\Phi(U\times \rho\cdot\Gamma) \subset \D.
\]

Pick $r_2$ in the interval $(1,r_1)$ such that $m(r_2\D) \subset U$.  Then for all $\la\in r_2\D$ and $z\in \rho\cdot\Gamma$,
\[
|\Phi(m(\la),z)| < 1.
\]
In particular, for all $\la\in r_2\D$,
\beq\label{Phi<1}
|\Phi(m(\la),f(\la))| < 1.
\eeq
Thus $\Phi\circ(m,f)$ belongs to the Schur class, and
\[
\Phi\circ(m,f)(\la_j) = \Phi\circ(m,h)(\la_j) =q(\la_j) \quad \mbox{ for } j=1,\dots,k.
\]
Hence $\Phi\circ(m,f)$ is a solution of the solvable Nevanlinna-Pick problem
\[
\la_j \mapsto q(\la_j), \quad j =1,\dots, k,
\]
as is $q\in \B l_{k-1}$.  Any $k$-point Nevanlinna-Pick problem that is solved by an element of $\B l_{k-1}$ is  extremally solvable and has a unique solution, and so $\Phi\circ(m,f) = q$.  This is a contradiction, since, by inequality \eqref{Phi<1}, $\Phi\circ(m,f)$ maps $r_2\D$ into $\D$, whereas $q$ maps $r_2\D\setminus \Delta$ to the complement of $\Delta$.   Hence $h$ is $k$-extremal. 
\end{proof}

If  Conjecture \ref{conj3} is true then all $n$-extremals for $\Gamma$ lie in $\E_{n-2, n}$.

\begin{observation} \label{C-n-extrem} 
Let $n \ge 2$. If condition $\mathcal{C}_{n -2}$ suffices for the solvability of $n$-point $\Gamma$-interpolation problems then  every rational $\Gamma$-inner function $h$ which is $n$-extremal for $\hol(\D,\Gamma)$ belongs to $\E_{n-2, n}$.
\end{observation}
\begin{proof} Let $h$ be $n$-extremal for $\hol(\D,\Gamma)$ and suppose that $h \notin \E_{n-2, n}$. Thus for every $\up \in \B l_{n-2}$ the function $\Phi \circ (\up,h) \notin \B l_{n-1}.$

Consider $n$ distinct points $\la_1, \dots, \la_n$ in $\D$ and let $z_j = h(\la_j)$, $1 \le j \le n$. Then, for all $\up \in \B l_{n-2}$, the Nevanlinna-Pick data
\begin{equation}
\label{up-ext-unsolv}
\la_j \mapsto \Phi(\up(\la_j), z_j), \quad 1 \le j \le n,
\end{equation}
are not extremally solvable, so that $\|X(\up) \| < 1$, where
the operator $X(\up)$ is defined on 
\begin{equation}\label{defM_1}
\mathcal{M} = \mathrm{span}\ \{K_{\lambda_1}, \dots,
K_{\lambda_n}\} 
\subset H^2,
\end{equation} 
by
\beq\label{defX_1}
X(\up)K_{\la_j}= \overline{ \Phi(\up(\la_j),z_j)}K_{\la_j}, \qquad 1\leq j\leq n.
\eeq
Thus, by the compactness of  $\B l_{n-2}$ and the continuity of $X(\cdot)$, there exists a positive constant $c < 1$ such that
\begin{equation}
\label{Xu_1-le1}
  \sup_{\up \in \B l_{n-2}} \|X(\up) \|\le c < 1.
\end{equation}

For $r \ge 1$ define the operator
$X_r(\up)$ on 
\begin{equation}\label{defM_r}
\mathcal{M}_r \df \mathrm{span}\ \{K_{\lambda_1/r}, \dots,
K_{\lambda_n/r}\} 
\subset H^2,
\end{equation} 
by
\beq\label{defX_r}
X_r(\up)K_{\la_j/r}= \overline{ \Phi(\up(\la_j/r),z_j)}K_{\la_j/r}, \qquad 1\leq j\leq n.
\eeq

\begin{lemma}\label{normX_r-less-1} There exists $r >1$ such that 
\[
\sup_{\up \in \B l_{n-2}} \|X_r(\up) \| < 1.
\]
\end{lemma}
\begin{proof} Note that
\begin{equation}\label{X_rtoT_r}
X_r(\up) = T_r Y_r(\up) T_r^{-1}
\end{equation} 
where 
\beq\label{defT_r}
T_r: \mathcal{M} \to \mathcal{M}_r: K_{\la_j} \mapsto K_{\la_j/r}, \qquad 1\leq j\leq n.
\eeq
and 
\beq\label{defY_r}
Y_r(\up): \mathcal{M} \to \mathcal{M}: K_{\la_j} \mapsto
\overline{ \Phi(\up(\la_j/r),z_j)}K_{\la_j}, \qquad 1\leq j\leq n.
\eeq
Since $Y_r(\up)$ is a diagonal operator on $ \mathcal{M}$ (with respect to a fixed basis) and its diagonal entries are continuous in $ r \ge 1$, $\up \in \B l_{n-2}$, the map $(r, \up) \mapsto Y_r(\up)$ is a continuous $\mathcal{L}(\mathcal{M})$-valued map on $[1, \infty) \times \B l_{n-2}$. 
Note that $Y_1(\up) = X(\up)$ and $\|X(\up) \|\le c < 1$ for all
$\up \in \B l_{n-2}$.
Thus there is $r_0 > 1$ such that 
\[
\|Y_r(\up) \| \le \tfrac{1}{2} (c+1) < 1
\]
for all $\up \in \B l_{n-2}$ and $1 < r <r_0$. Therefore
\[
\|X_r(\up) \|  = \|T_r Y_r(\up) T_r^{-1} \| 
\le \tfrac{1}{2} (c+1)\|T_r \|\|T_r^{-1} \| 
\]
for all $\up \in \B l_{n-2}$ and $1 < r <r_0$.
Since $K_{\la_j/r} \to K_{\la_j}$ in $ H^2$ as $r \to 1$, it 
is straightforward to show that 
 $\|T_r \| \to 1$ and $\|T_r^{-1} \| \to 1$ as $r \to 1$. Hence, for $r >1$ sufficiently close to $1$, we have $\|X_r(\up)\| < 1$ for all $\up \in  \B l_{n-2}$.\end{proof}

{\bf Conclusion of the Proof of Observation \ref{C-n-extrem}.} By Lemma \ref{normX_r-less-1}, there exists  $r >1$ such that
$\|X_r(\up) \|< 1$ for all $\up \in \B l_{n-2}$.
Therefore, for every 
Blaschke product $\up \in \B l_{n-2}$, the Nevanlinna-Pick data
\begin{equation}
\label{up-ext-unsolv_r}
\la_j/r \mapsto \Phi(\up(\la_j/r), z_j), \quad 1 \le j \le n,
\end{equation}
are solvable.
In other words, $\la_j/r \mapsto z_j$, $ 1 \le j \le n$, satisfy the condition $\mathcal{C}_{n-2}$.

By assumption, $\mathcal{C}_{n -2}$ suffices for solvability of $n$-point $\Gamma$-interpolation problems. Therefore there exists $f \in \hol(\D, \Gamma)$ such that 
\[
f(\la_j/r) = z_j, \quad 1 \le j \le n.
\]
Then the function $g(\la)= f (\la/r)$ belongs to $\hol(r\D, \Gamma)$ and satisfies
\[
g(\la_j) = f(\la_j/r) = z_j= h(\la_j), \quad 1 \le j \le n.
\]
This contradicts the $n$-extremality of $h$. Thus $h \in \E_{n-2, n}$.
\end{proof}

\section{Complex geodesics of $\G $ and the classes $\E_{\nu 2}$}\label{geodesics}

In this section we shall show that all the classes $\mathcal{E}_{\nu 2}$  consist of the superficial rational $\Gamma$-inner functions together with the complex geodesics of $\G$.
First we recall a result from \cite{AY06}.
\begin{proposition}\label{complex-geodesics-Phi}   An analytic  function $h:\D \to \G  $ is  a  complex geodesic of $\G$ if and only if there is an $\omega \in \T$ such that $\Phi_{\omega} \circ h \in\Aut \mathbb{D}$.
Furthermore, every complex geodesic of $\G$ is $\Gamma$-inner.
\end{proposition}
\begin{proof} $\Leftarrow$ Suppose there is $\omega \in \T$ such that $\Phi_{\omega} \circ h \in \Aut \mathbb{D}$, say $\Phi_{\omega} \circ h = \up \in \Aut \mathbb{D}$. Then $g=\up^{-1}  \circ \Phi_{\omega}:  \G \to \D$ is an analytic left inverse of $h$. Therefore, $h$ is  a  complex geodesic of $\G$.

$\Rightarrow$ Let $h$ be a  complex geodesic of $\G$ and let $g : \G \to \D$ be
an analytic left inverse of $h$. For any distinct points $\lambda_1, \lambda_2 \in \D$, we have
\[
C_{\G}(h(\lambda_1), h(\lambda_2)) \le \rho(\lambda_1, \lambda_2).
\]
On the other hand, since $ g \circ h = \id_{\D}$,
\[
\rho(\lambda_1, \lambda_2) = \rho(g \circ h(\lambda_1), g \circ h(\lambda_2)) \le C_{\G}(h(\lambda_1), h(\lambda_2)). 
\]
Therefore, 
\[
C_{\G}(h(\lambda_1), h(\lambda_2)) = \rho(\lambda_1, \lambda_2).
\]
By \cite[Theorem 1.2]{AY06}, there is there is $\omega \in \T$ such that $\Phi_{\omega} \circ h \in\Aut \mathbb{D}$.

By \cite[Lemma 1.1]{AY06}, $h$ is $\Gamma$-inner.
\end{proof}

\begin{corollary}\label{E_0_2=complex-geodesics}
The set $\E_{0 2}\setminus \E_{0 1}$ is precisely the set of complex geodesics of $\G$. 
\end{corollary}
\begin{proof} 
Let $h \in \E_{0 2}\setminus \E_{0 1}$.
By Definition \ref{E_nu_k},  there exists $\omega \in \T$ such that 
\[
\Phi_{\omega}\circ h= \frac{2 \omega p -s}{2- \omega s}  \in \B l_{1}\setminus \B l_0 = \Aut \D.
\]
Hence, by Proposition \ref{complex-geodesics-Phi}, $h$ is a complex geodesic of $\G$.

Conversely, suppose that $h$ is a complex geodesic of $\G$.  By Proposition \ref{complex-geodesics-Phi}, there exists $\omega\in \T$ such that $\Phi_\omega \circ h \in\Aut\D$, and hence $h\in \E_{02}$.  Since  $h(\D) \subset \G$, $h$ is not superficial, and so $h \notin \E_{01}$.
\end{proof}

\begin{theorem}\label{Enu2}
For $\nu\geq 0$ the set $\mathcal{E}_{\nu 2}$ is the union of the set of superficial rational $\Gamma$-inner functions and the set of complex geodesics of $\G$.
\end{theorem}
\begin{proof}
Since $\mathcal{E}_{0 2}\subset\mathcal{E}_{\nu 2}$, it follows from Corollary \ref{E_0_2=complex-geodesics} that $\mathcal{E}_{\nu 2}$ contains all superficial  rational $\Gamma$-inner functions and all complex geodesics of $\G$.  If $h\in \mathcal{E}_{\nu 2}$ then, by Lemma \ref{notsuper}, either $h$ is superficial or $h\in\hol(\D,\G)$.  In the latter case, $h$ is $2$-extremal by Theorem \ref{Enukextrem} and $h$ is a complex geodesic of $\G$ by Corollary \ref{compgeos}. 
\end{proof}

\section{Condition $\mathcal{C}_\nu$ and the classes $\E_{\nu k}$}\label{connections-C-E}

It is clear that $\mathcal{C}_{\nu}(\la,z)$ implies $\mathcal{C}_{\nu-1}(\la,z)$ for any $\Gamma$-interpolation data $\la\mapsto z$. 
 To show that $\mathcal{C}_\nu$ is {\em strictly stronger} than $\mathcal{C}_{\nu-1}$ we need to find data 
\begin{equation}
\label{GdataEC}
\la=(\la_1,\dots, \la_k), \quad z=(z_1, \dots, z_k),
\end{equation}
where $\lambda_1, \dots, \lambda_k$ are distinct points in $\D$
and $z_j = (s_j,p_j) \in \G$ for $j=1,\dots,k$,
such that 

(i) for every Blaschke product $\up$ of degree at most $\nu-1$,
\begin{equation}\label{upsdata-nu-1}
\la_j \mapsto \frac{2\up(\la_j)p_j -s_j}{2-\up(\la_j)s_j}, \quad j=1,\dots,k,
\end{equation}
are solvable Nevanlinna-Pick data, but 

(ii) there is a Blaschke product $m$ of degree $\nu$ such that
\begin{equation}
\label{upsdata-nu}
\la_j \mapsto \frac{2 m(\la_j)p_j -s_j}{2-m(\la_j)s_j}, \quad j=1,\dots,k,
\end{equation}
are not solvable Nevanlinna-Pick data.

For distinct points $\la_1, \dots,\la_k$  in $\D$, we define
\[
\mathrm {Solv}(\la_1, \dots,\la_k) = \{(f(\la_1), \dots, f(\la_k)) \in \D^k: f \in \Schur\},
\]
and 
\[
\mathrm {Unsolv}(\la_1, \dots,\la_k) = \C^k \setminus {\rm Solv}(\la_1, \dots,\la_k).
\]
Thus $w =(w_1, \dots, w_k) \in {\rm Solv}(\la_1, \dots,\la_k)$ if and only if $ \la_j \mapsto w_j, \; j=1,\dots,k, $ are solvable Nevanlinna-Pick data. 
\begin{proposition}\label{Solv-closed}  Let   $\la_1, \dots,\la_n$  be distinct points in $\D$.

{\rm (i)} ${\rm Solv}(\la_1, \dots,\la_n)$ is closed in $\C^n$.

{\rm (ii)} Let $w=(w_1, \dots, w_n) \in {\rm Solv}(\la_1, \dots,\la_n)$. The  Nevanlinna-Pick data $\la_j \mapsto w_j$, $ j=1,\dots,n$, are  extremally solvable if and only if  $w \in \partial{\rm Solv}(\la_1, \dots,\la_n)$.
\end{proposition}

\begin{proof} Part (i) is immediate from Pick's Theorem, which asserts that
 $w=(w_1, \dots, w_n) \in {\rm Solv}(\la_1, \dots,\la_n)$ if and only if 
\[
\left[\frac{1-\bar{w_i} w_j}{1-\bar{\la_i} \la_j} \right]_{i,j=1}^n \ge 0.
\]

(ii) Suppose that the  Nevanlinna-Pick data $\la_j \mapsto w_j$, $ j=1,\dots,n$, are  extremally solvable. We will show that $w \in \partial{\rm Solv}(\la_1, \dots,\la_n)$ by induction on $n$.
It is true when $n=1$, since $\la_1 \mapsto w_1$ is an extremally 
solvable Nevanlinna-Pick datum if and only if $|w_1|=1$. Consider $n \ge 2$ and suppose (ii) holds for $n-1$. Let $\la_j \mapsto w_j$, $1 \le j \le n$, be  extremally solvable Nevanlinna-Pick data, and let 
 $f \in \Schur$ be a solution. The Schur reduction $f_1$ of $f$ at $\la_1$,
\begin{equation} \label{SchRed}
f_1(\la) = \frac{1-\bar{\la}_1 \la}{\la - \la_1} \frac{f(\la)- w_1}{1 - \bar{w}_1 f(\la)},
\end{equation} 
also lies in $\Schur$ and satisfies $f_1(\la_j) = w'_j$, $2 \le j \le n$,
where
\[
 w'_j =\frac{1-\bar{\la}_1 \la_j}{\la_j - \la_1} \frac{w_j- w_1}{1 - \bar{w}_1 w_j}.
\]
We claim that the Nevanlinna-Pick data 
\begin{equation} \label{ext-solv-NP}
\la_j \mapsto w'_j, \quad 2 \le j \le n,
\end{equation}
are also extremally solvable. They are certainly solvable, since $f_1$ is a solution. If they are not  extremally solvable then there are two distinct functions $f_1, \tilde{f}_1 \in \Schur$ that solve the data, and on inverting the relation \eqref{SchRed} we obtain two distinct solutions of the Nevanlinna-Pick problem with data $\la_j \mapsto w_j$, $ 1 \le j \le n$;
this contradicts the fact that extremally solvable Nevanlinna-Pick problems have unique solutions. The claim follows.

By the inductive hypothesis there is a sequence $(w'(k))_{k \ge 1}$ in
${\rm Unsolv}(\la_2, \dots,\la_n)$ converging to $w'= (w'_2, \dots, w'_n)$ as $k \to \infty$. Let
\[
 w_j(k) =\frac{\frac{\la_j - \la_1}{1-\bar{\la}_1 \la_j} w'_j(k) +w_1}
{1 + \bar{w}_1 \frac{\la_j - \la_1}{1-\bar{\la}_1 \la_j}w'_j(k)}
\;\; \text{for} \;\; k \ge 1,\; 2 \le j \le n.
\]
Then 
\[
w(k) \df (w_1, w_2(k), \dots, w_n(k)),\;\; k =1,2,\dots,
\]
is a sequence in ${\rm Unsolv}(\la_1, \la_2, \dots,\la_n)$ that converges to $w$ as $k \to \infty$. Thus $w \in \partial{\rm Solv}(\la_1, \la_2, \dots,\la_n)$.

Hence, by induction, $w \in \partial{\rm Solv}(\la_1, \dots,\la_n)$  for all $n \in \N$.

Suppose the  Nevanlinna-Pick data $\la_j \mapsto w_j$, $ j=1,\dots,n$, are not extremally solvable. Then
\[
\left[\frac{1-\bar{w_i} w_j}{1-\bar{\la_i} \la_j} \right]_{i,j=1}^n > 0.
\]
Note that $|w_j| <1$ for all $ j=1,\dots,n$. Hence, by continuity, there are 
neighbourhoods $U_1,  \dots, U_n$ of $w_1,  \dots, w_n$ in $\D$ such that, for all $z_1 \in U_1$, $\dots$, $z_n \in U_n$, we have
\[
\left[\frac{1-\bar{z_i} z_j}{1-\bar{\la_i} \la_j} \right]_{i,j=1}^n > 0.
\]
Hence  $(z_1, z_2, \dots, z_n)  \in {\rm Solv}(\la_1, \dots,\la_n)$ for all 
$(z_1, z_2, \dots, z_n) \in U_1 \times \dots \times U_n$. That is, there is a 
neighbourhood $U_1 \times \dots \times U_n$ of $(w_1,  \dots, w_n)$ that is contained in ${\rm Solv}(\la_1, \dots,\la_n)$. In other words   $(w_1,  \dots, w_n)$ is an interior point of ${\rm Solv}(\la_1, \dots,\la_n)$, hence is not in $\partial{\rm Solv}(\la_1, \dots,\la_n)$.
\end{proof}

\begin{proposition}\label{E_nuandC_nu}
If there exists a nonconstant  function $h \in \E_{\nu k} \setminus \E_{\nu-1, k}$ then $\mathcal{C}_\nu$ is strictly stronger than $\mathcal{C}_{\nu-1}$. In fact there is a set of  $\Gamma$-interpolation data $\lambda_j \mapsto z_j$ with $k$ interpolation points which satisfies $\mathcal{C}_{\nu -1}$ but not $\mathcal{C}_{\nu}$.
\end{proposition}
It follows of course from Theorem \ref{necG} that the data $\lambda_j \mapsto z_j$ are not solvable.

\begin{proof} Pick any $k$  distinct points $\la_1, \dots,\la_k$ in $\D$ and let $h(\la_j) = (s_j, p_j)$, $ j=1,\dots,k$.
Since $h=(s,p) \notin \E_{\nu-1, k}$, there is no 
$\up \in \B l_{\nu-1}$ such that 
\begin{equation} \label{not-in-Bl(k-1)}
\frac{2 \up p -s}{2- \up s}  \in \B l_{k-1}.
\end{equation}
Consider any $\up \in \B l_{\nu-1}$. The Nevanlinna-Pick data 
\begin{equation} \label{solv-NP}
\la_j \mapsto \Phi(\up(\la_j),s_j,p_j)=\frac{2\up(\la_j)p_j -s_j}{2-\up(\la_j)s_j}, \quad j=1,\dots,k,
\end{equation}
are solvable, since $\la \mapsto \Phi(\up(\la), h(\la))$ is in $\Schur$ and satisfies the interpolation condition \eqref{solv-NP}. However, if the data \eqref{solv-NP} are extremally solvable, then the function 
$\Phi \circ (\up, h)$ is a Blaschke product of degree at most $k-1$ \cite[Theorem 6.15]{AgMcC}, contrary to the hypothesis that there is no 
$\up \in \B l_{\nu-1}$ such that equation $ \Phi \circ (\up, h) \in \B l_{k-1}$.
Therefore, for all $\up \in \B l_{\nu-1}$, the data \eqref{solv-NP}
are not extremally solvable and so $(w_1, \dots, w_k)$ lies in the interior of ${\rm Solv}(\la_1, \dots,\la_k)$, where $w_j=\Phi(\up(\la_j),s_j,p_j)$. Recall from Section \ref{extremality} that, by Pick's Theorem, the Nevanlinna-Pick data
\begin{equation}\label{sol-not-extr}
\la_j \mapsto \Phi(\up(\la_j),s_j,p_j), \qquad 1 \leq j \leq k,
\end{equation}
are solvable if and only if the operator
$X(\up)$ on $\mathcal{M}$ given by
\beq\label{CEdefX}
X(\up) K_{\la_j}= \overline{\Phi(\up(\la_j),s_j,p_j)} K_{\la_j}, \qquad 1\leq j\leq k,
\eeq
is a contraction. The Nevanlinna-Pick data \eqref{sol-not-extr} are extremally solvable if and only if
$$
\sup_\upsilon \|X(\upsilon)\| = 1
$$
where the supremum is over all Blaschke products of degree at most $\nu-1$.
Thus, for all $\up \in \B l_{\nu -1}$,
$$
   \|X(\up) \| < 1.
$$
By the compactness of $\B l_{\nu-1}$ and the continuity of $X$, there is a positive constant $c <1$
such that, for all $\up \in \B l_{\nu -1}$,
$$
   \|X(\up) \| \le c < 1.
$$
Hence there is a neighbourhood  $U$ of $(s_1, \dots, s_k)$ in $(2\D)^k$
such that, for all 
\newline $(\tilde{s}_1, \dots, \tilde{s}_k) \in U$ and all $\up \in \B l_{\nu -1}$,
\begin{equation}\label{sol-in-U}
\la_j \mapsto \Phi( \up(\la_j),\tilde{s}_j,p_j )=
\frac{ 2\up(\la_j)p_j -\tilde{s}_j }{ 2-\up(\la_j) \tilde{s}_j }, \quad j=1,\dots,k,
\end{equation}
are solvable  Nevanlinna-Pick data.

By assumption,  $h=(s,p) \in \E_{\nu k}$, and so there is 
$ m \in \B l_{\nu}$ such that 
\[
q \df \Phi \circ (m,h) = \frac{2 m p -s}{2- m s}  \in \B l_{k-1}.
\]
Let 
\[
\tilde{w}_j \df q(\la_j) = \Phi( m(\la_j),s_j,p_j )=
\frac{2 m(\la_j) p_j - s_j}{2- m(\la_j)s_j}, \quad j=1,\dots,k.
\]
Then the Nevanlinna-Pick data 
\[
\la_j \mapsto \tilde{w}_j, \quad j=1,\dots,k,
\]
are solvable and $q$ is a solution. By \cite[Lemma 6.19]{AgMcC}, the Pick matrix
\[
\left[\frac{1-\overline{\tilde{w}_i} \tilde{w}_j}{1-\bar{\la_i} \la_j} \right]_{i,j=1}^k \ge 0.
\]
is positive and of rank at most $ k-1$. Hence the Pick matrix is singular and the data $\la_j \mapsto \tilde{w}_j$, $j=1,\dots,k,$
are extremally solvable. By Proposition \ref{Solv-closed}, 
\[
\tilde{w} =  \left( \Phi( m(\la_j),s_j,p_j ) \right)_{j=1}^k \in \partial{\rm Solv}(\la_1, \dots,\la_k).
\]
Define an analytic function
\[
F: U \to  \D^k: (\zeta_1,\dots, \zeta_k) \mapsto \left( \Phi( m(\la_j), \zeta_j,p_j ) \right)_{j=1}^k.
\]
Note that $F(s_1, \dots, s_k) = \tilde{w} \in \partial{\rm Solv}(\la_1,  \dots,\la_k)$. The Jacobian matrix
\[
J_F(s_1, \dots, s_k) = \diag\left( 2 \frac{ m(\la_j)^2 p_j - 1}{(2- m(\la_j)s_j)^2}\right)_{j=1}^k
\]
is nonsingular. Hence, by the Inverse Function Theorem \cite{Dieu}, there is a neighbourhood  $W$ of $(s_1, \dots, s_k)$ in $U$
such that $F(W)$ is an open neighbourhood of $\tilde{w}$ in $\C^k$ and  $F|_W$ is bijective.
Pick a point $w' \in F(W) \setminus {\rm Solv}(\la_1, \dots,\la_k)$ and let $(\tilde{s_1},\dots, \tilde{s_k})= F^{-1}(w')$.
Then $(\tilde{s_1},\dots, \tilde{s_k}) \in W \subset U$, and so 
\begin{equation}\label{proof-upsdata-nu-1}
\la_j \mapsto \frac{2\up(\la_j)p_j -\tilde{s_j}}{2-\up(\la_j)\tilde{s_j}}, \quad j=1,\dots,k,
\end{equation}
are solvable Nevanlinna-Pick data for  every Blaschke product $\up \in \B l_{\nu-1}$. Thus $(\la_j, \tilde{s_j}, p_j)$, $j=1,\dots,k$, satisfy $\mathcal{C}_{\nu-1}$.

On the other hand, since $F(\tilde{s_1},\dots, \tilde{s_k})= w' \in {\rm Unsolv}(\la_1, \dots,\la_k)$,
\begin{equation}
\label{proof-upsdata-nu}
\la_j \mapsto \frac{2 m(\la_j)p_j -\tilde{s_j}}{2-m(\la_j)\tilde{s_j}}, \quad j=1,\dots,k,
\end{equation}
are not solvable Nevanlinna-Pick data for $m \in \B l_{\nu}$, and thus   $(\la_j, \tilde{s_j}, p_j)$, $j=1,\dots,k$, do not satisfy $\mathcal{C}_{\nu}$.
\end{proof}

\section{Inequations for the classes $\E_{\nu k}$}\label{inequations-for-E}

In order to apply Proposition \ref{E_nuandC_nu} we must establish the strict inclusion 
\[
\mathcal{E}_{\nu -1,k} \subsetneq \mathcal{E}_{\nu,k}
\]
for a suitable $k$.

\begin{proposition}\label{E_nu-munis-E_nu-1}  For all $\nu \ge 1$ and $0<r<1$,  the function 
\begin{equation}\label{h_nu}
h_{\nu}(\la) =\left( 2(1-r) \frac{\la^{\nu+1}}{1+ r\la^{2\nu+1}},
\frac{\la(\la^{2\nu+1}+r)}{1+ r\la^{2\nu+1}} \right), \;\la \in \D,
\end{equation}
belongs to $\E_{\nu, \nu +2} \setminus \E_{\nu-1, \nu +2}$.
\end{proposition}

We require two lemmas.

\begin{lemma}\label{h-Ginner} Let $h$ be analytic on $\Delta$ and let $h(\T) \subset b \Gamma$. Then $h(\Delta) \subset \Gamma$ and so $h$ is $\Gamma$-inner.
\end{lemma}
\begin{proof}   Let $h=(s,p)$.  Observe that $|s| \le 2$ on $\Delta$, by the Maximum Principle, and $p$ is inner. We can suppose that $s^2 -4p$ is not identically $0$. 

Let the zeros of $s^2 -4p$ on $\T$ be $\la_1, \dots, \la_N$ and let $h(\la_j)= (2 \bar{\omega_j}, \bar{\omega_j}^2)$, $\omega_j \in \T$. Consider any $\omega \in \T \setminus \{\omega_1, \dots, \omega_n \}$.
Note that $\Phi_{\omega}$ is analytic on a neighbourhood of $\Gamma \setminus \{2\bar{\omega}, \bar{\omega}^2 \}$.

For every $\la \in  \T \setminus \{\la_1, \dots, \la_N\}$, we have $|s(\la)| < 2$. If $|s(\la)| = 2$ for $\la \in  \T$ then $h(\la) \in \{(2 \bar{\omega_j}, \bar{\omega_j}^2): \ 1 \le j \le N\}$. Hence $2-\omega s(\la) \neq 0$ for all $\la \in  \T$. Therefore
\[
\Phi_\omega \circ h =\frac{ 2 \omega p -s}{2 - \omega s}
\]
is analytic on $\Delta$. 

For any $\la \in  \T$, $h(\la) \in b \Gamma$ and hence $|\Phi_\omega \circ h (\la)|=1$. By the Maximum Principle, $|\Phi_\omega \circ h (\la)|\le 1$
on $\Delta$.  Since this is true for all but finitely many $\omega\in\T$, by Proposition \ref{critGam}(2), $h(\la)\in\Gamma$.
\end{proof}

\begin{lemma}\label{canc-degree-up} 
Let $h =(s,p)$ be a rational $\Gamma$-inner function. Suppose that\\
{\rm (i)} $h$ has $N$ distinct royal nodes $\omega_j$, $1\le j \le N$, on $\T$, and\\
{\rm (ii)} there is a finite Blaschke product $m$ of degree at most $\frac{1}{2}N$ such that $\Phi \circ (m, h)$ has cancellations at $\omega_j$, $1\le j \le N$.\\
Then, for any  $\up \in \B l_{d(m)-1}$, the function 
$\Phi \circ (\up, h)$ has no more than $d(m) +d(\up)$  cancellations  and the degree of $\Phi \circ (\up, h)$ is at least $d(p) - d(m)$.
\end{lemma}

\begin{proof}    
By Proposition \ref{canc=royal}, since there are cancellations in $\Phi \circ (m, h)$ at $\omega_j$, $1\le j \le N$, we have  $m(\omega_j)=\frac{1}{2}\overline{ s(\omega_j)}$, $1\le j \le N$.

Suppose that, for some $\up \in \B l_{d(m)-1}$, the function 
$\Phi \circ (\up, h)$ has  $d(m) +d(\up)+1$  cancellations. Then 
$\up(\omega_j)=\frac{1}{2}\overline{ s(\omega_j)}$ at $d(m) +d(\up)+1$ distinct points $\omega_j$ on $\T$.
Hence the rational function $m- \up$ vanishes at $d(m) +d(\up)+1$ distinct points $\omega_j$ on $\T$. The degree of $m- \up$ is at most $ d(m) +d(\up)$ and so $m=\up$. This is a contradiction to the assumption that $\up \in \B l_{d(m)-1}$.

By Theorem \ref{crit-canc=royal}, $\Phi \circ (\up,h)$ has no double cancellations at $\omega_j$. Thus $\Phi \circ (\up,h)$ has degree at least $ d(\up p) - (d(m)+ d(\up))=d(p) - d(m)$.
\end{proof}

We can now prove Proposition \ref{E_nu-munis-E_nu-1}.

\begin{proof}  It is clear that $h_\nu$ is analytic on $\Delta$.   Let $h_\nu=(s,p)$.
It is simple to check that $s=\bar{s}p$ on $\T$, 
 that $|s| \le 2$ on $\T$ and that
$|s(\la)| = 2$ if and only if $\la^{2\nu+1}= -1$.
For all $\la \in \T$,
\[
|p(\la)| =\left| \frac{\la(\la^{2\nu+1}+r)}{1+ r\la^{2\nu+1}} \cdot\frac{1}{\bar{\la}^{2\nu+1}} \right|=\left| \frac{\la^{2\nu+1}+r}{\bar{\la}^{2\nu+1}+ r}  \right| =1.
\]
By Proposition \ref{critGam}, we have $h_\nu(\T) \subset b\Gamma$. Thus,
by Lemma \ref{h-Ginner}, $h_{\nu}$ is $\Gamma$-inner.

Let $m(\la) = - \la^\nu$, so that $m\in \B l_\nu$.
It is simple to verify that
\[
\Phi \circ (m,h_\nu)=\frac{2 m p -s}{2- m s}(\la)=- \la^{\nu+1} \in \B l_{\nu+1},
\]
and so $h_\nu \in \E_{\nu, \nu +2}$.

To prove that  $h_\nu$ is not in  $ \E_{\nu-1, \nu +2}$ we must show that, for all $\up \in \B l_{\nu-1}$, the Blaschke product $\Phi \circ (\up,h_\nu)$  has degree at least $ \nu +2$.
By Proposition \ref{canc=royal}, for $\up \in \B l_{\nu-1}$, if the function
\[
\Phi \circ (\up, h_\nu) = \frac{2 \up p -s}{2- \up s} 
\]
has a cancellation at $\zeta$, then $\zeta \in \T$, $\zeta$ is a royal node for $h_\nu$ and $|s(\zeta)|=2$. 
The royal nodes for $h_\nu$, being the points at which $|s|=2$, are the $(2\nu +1)$th roots of $-1$, that is, 
\[
\omega_j = e^{\ii \pi (2j+1)/(2\nu+1)}, \; j= 0, \dots, 2 \nu.
\]
Note that $s(\omega_j) =2 \omega_j^{\nu +1}$, $j= 0, \dots, 2 \nu$.
For the finite Blaschke product $m(\la) = - \la^\nu$, we have
$m(\omega_j)= -\omega_j^\nu=\overline{\omega_j}^{\nu +1}$ since $\omega_j^{2\nu +1}=-1$. Hence $m(\omega_j)= \frac{1}{2}\overline{s(\omega_j)}$,
 $j= 0, \dots, 2 \nu$.
Thus $\Phi \circ (m, h)$ has cancellations at  $\omega_j$, $0\le j \le 2\nu$.
By Lemma \ref{canc-degree-up}, for every  $\up \in \B l_{\nu-1}$, the function $\Phi \circ (\up,h_\nu)$ has no more than $d(\up)+ \nu$ different points of cancellation and the function 
 $\Phi \circ (\up,h_\nu)$ has degree at least $ d(p)- d(m)=  2\nu +2 -\nu=  \nu + 2$. Therefore $h_\nu$ is not in  $ \E_{\nu-1, \nu +2}$.
\end{proof}

Our main theorem follows easily.

\begin{theorem}\label{C_nu-stronger-C_nu-1} For all $\nu \ge 1$,
the condition $\mathcal{C}_\nu$ is strictly stronger than $\mathcal{C}_{\nu-1}$. In fact there is a set of $\Gamma$-interpolation data $\lambda_j \mapsto z_j$ with $\nu +2$ interpolation points which satisfies $\mathcal{C}_{\nu -1}$ but not $\mathcal{C}_{\nu}$.
\end{theorem}

\begin{proof} By Proposition \ref{E_nu-munis-E_nu-1}, there exists a nonconstant  function $h \in \E_{\nu, \nu +2} \setminus \E_{\nu-1, \nu +2}$. By Proposition \ref{E_nuandC_nu},
the condition $\mathcal{C}_\nu$ is strictly stronger than $\mathcal{C}_{\nu-1}$, and furthermore, there is a set of  $\Gamma$-interpolation data $\lambda_j \mapsto z_j$ with $\nu +2$ interpolation points which satisfies $\mathcal{C}_{\nu -1}$ but not $\mathcal{C}_{\nu}$.
\end{proof}

As we observed above, $\mathcal{C}_0$ is necessary and sufficient for solvability of a $\Gamma$-interpolation problem when $n=2$, but a consequence of Theorem \ref{C_nu-stronger-C_nu-1} is:

\begin{corollary}
For all $n \ge 3$,
Condition $\mathcal{C}_{n-3}$ does not suffice for the solvability of an $n$-point $\Gamma$-interpolation problem.
\end{corollary}

\section{Table of relations between the classes $\E_{\nu k}$}\label{summary-E-cl}

The following table summarises the relations between $\E$-classes established above.

{\small
\begin{equation*}
\begin{array}{cccccccccccccc}
\E_{0 1} &   \stackrel {\rm (4,5)} 
{\subsetneq} & \E_{0 2} &  \stackrel {\rm (1)} 
{\subsetneq} &  \E_{0 3} & \stackrel {\rm (1)} 
{\subsetneq}  & \E_{0 4}  &  \stackrel {\rm (1)} 
{\subsetneq}  & \E_{0 5}  & \stackrel {\rm (1)} 
{\subsetneq}  & \E_{0 6}  & \stackrel {\rm (1)} {\subsetneq }
 & \E_{0 7}  & \subsetneq  \dots\\
\stackrel {\rm (4)} {\|} & ~ & \stackrel {\rm (5)} {\|}  & ~ &  \stackrel {\rm (2)} {\not{|\bigcap}} & ~ & \bigcap & ~ &  \bigcap
& ~ &\stackrel {\rm (3)} {\not{|\bigcap}}& ~ &\stackrel {\rm (3)} {\not{|\bigcap}}& ~ \\
\E_{1 1} &   \stackrel {\rm (4,5)} 
{\subsetneq} & \E_{1 2} &  \stackrel {\rm (1)} 
{\subsetneq} &  \E_{1 3} & \stackrel {\rm (1)} 
{\subsetneq}  & \E_{1 4}  &  \stackrel {\rm (1)} 
{\subsetneq}  & \E_{1 5}  & \stackrel {\rm (1)} 
{\subsetneq}  & \E_{1 6}  & \stackrel {\rm (1)} {\subsetneq }
 & \E_{1 7}  & \subsetneq  \dots\\
\stackrel {\rm (4)} {\|}  & ~ & \stackrel {\rm (5)} {\|}  & ~ &  \bigcap & ~ & \stackrel {\rm (2)} {\not{|\bigcap}} & ~ &  \bigcap & ~ &  \bigcap & ~ &  \bigcap \\
\E_{2 1} &   \stackrel {\rm (4,5)} 
{\subsetneq} & \E_{2 2} &  \stackrel {\rm (1)} 
{\subsetneq} &  \E_{2 3} & \stackrel {\rm (1)} 
{\subsetneq}  & \E_{2 4}  &  \stackrel {\rm (1)} 
{\subsetneq}  & \E_{2 5}  & \stackrel {\rm (1)} 
{\subsetneq}  & \E_{2 6}  & \stackrel {\rm (1)} {\subsetneq }
 & \E_{2 7}  & \subsetneq  \dots\\
\stackrel {\rm (4)} {\|}  & ~ & \stackrel {\rm (5)} {\|}  & ~ &  \bigcap & ~ & \bigcap & ~ & \stackrel {\rm (2)} {\not{|\bigcap}} & ~ &  \bigcap & ~ &  \bigcap \\
\E_{3 1} &   \stackrel {\rm (4,5)} 
{\subsetneq} & \E_{3 2} &   
{\subset} &  \E_{3 3} & \stackrel {\rm (1)} 
{\subsetneq}  & \E_{3 4}  &  \stackrel {\rm (1)} 
{\subsetneq}  & \E_{3 5}  & \stackrel {\rm (1)} 
{\subsetneq}  & \E_{3 6}  & \stackrel {\rm (1)} {\subsetneq }
 & \E_{3 7}  & \subsetneq  \dots\\
\stackrel {\rm (4)} {\|}  & ~ & \stackrel {\rm (5)} {\|}  & ~ &  \bigcap & ~ & \bigcap & ~ & \bigcap & ~ &  \stackrel {\rm (2)} {\not{|\bigcap}}& ~ &  \bigcap \\
\E_{4 1} &   \stackrel {\rm (4,5)} 
{\subsetneq} & \E_{4 2} &  
{\subset} &  \E_{4 3} & 
{\subset}  & \E_{4 4}  &  \stackrel {\rm (1)} 
{\subsetneq}  & \E_{4 5}  & \stackrel {\rm (1)} 
{\subsetneq}  & \E_{4 6}  & \stackrel {\rm (1)} {\subsetneq }
 & \E_{4 7}  & \subsetneq  \dots\\
\stackrel {\rm (4)} {\|}  & ~ & \stackrel {\rm (5)} {\|}  & ~ &  \bigcap & ~ & \bigcap & ~ & \bigcap & ~ &  \bigcap & ~ & \stackrel {\rm (2)} {\not{|\bigcap}} \\
\dots &   ~ & \dots &   ~ & \dots &   ~ &\dots &   ~ &\dots &   ~ &\dots &   ~ &\dots &   ~ \\
 \E_{\nu 1} &   \stackrel {\rm (4,5)} 
{\subsetneq}  & \E_{\nu 2}  &  \subset & \E_{\nu 3} &  \subset & \E_{\nu 4}&  \subset &  \E_{\nu 5}&\subset & \E_{\nu 6}& \subset & \E_{\nu 7}& \subset \dots\\
\stackrel {\rm (4)} {\|}  & ~ & \stackrel {\rm (5)} {\|}  & ~ &  \bigcap & ~ & \bigcap & ~ & \bigcap & ~ &  \bigcap& ~ & \bigcap  \\
\dots &   ~ & \dots &   ~ & \dots &   ~ &\dots &   ~ &\dots &   ~ &\dots &   ~ &\dots &   ~ \\
\end{array}
\end{equation*}
}

\begin{remark} \label{Enotes}{\rm

(1) In Corollary \ref{h=(0,p)-E} we  proved that, for all $\nu \ge 0$, every  rational  $\Gamma$-inner function $h$ of the form $(0,p)$ is in $\E_{\nu, \nu+ d(p) +1} \setminus \E_{\nu, \nu + d(p)}$.
 In Example  \ref{E-1-k+2}, for a rational inner function  $\psi$ on $\D$, we considered the rational $\Gamma$-inner function
\[
h_\psi(\la) =\left( \la + \la \psi(\la),\la^{2} \psi(\la) \right), \;\la \in \D,
\]
 and proved that  $ h_\psi \in \E_{1, d(\psi)+2} \setminus \E_{1, d(\psi)+1}$.

(2) In Proposition \ref{E_nu-munis-E_nu-1}, for all $\nu \ge 1$, we presented  a  function $h_\nu \in \E_{\nu, \nu +2} \setminus \E_{\nu-1, \nu +2}$.

(3) In Example \ref{E-1-2j+4}, for $j =1,2, \dots$, we  constructed a rational  function
$h_j \in \E_{1, 2j+4} \setminus \E_{0,  2j+4}$.

(4) In Section \ref{superficial} we  showed that,  for every $\nu \ge 0$, the class $\E_{\nu 1}$  consists of the superficial rational $\Gamma$-inner functions.

(5) In Section \ref{geodesics} we  proved that,  for every $\nu \ge 0$, the class $\E_{\nu 2}$ comprises precisely the complex geodesics of $\G$ and the superficial rational $\Gamma$-inner functions.

(6) In Theorem \ref{Enukextrem}  we showed that, for $\nu \ge0$ and $k \ge 2$, functions in $\mathcal{E}_{\nu k}$  are either superficial or $k$-extremal.
}
\end{remark}

\begin{remark} \label{conj-extr}{\rm
If  the $\Gamma$-interploation Conjecture  is true then, for $k \ge 3$, the columns are all ultimately constant:
\[
\E_{k-2, k}  =  \E_{k-1, k}= \E_{k k}  =  \E_{k+1, k}= \dots.
\]
By Theorem \ref{Enukextrem},  if $h \in \E_{\nu, k}$ then either  $h$ is superficial, in which case $ h$ belongs to all $ \E_{\nu n}$, or $h$ is $k$-extremal, and then, by Observation \ref{C-n-extrem}, $h \in \E_{k-2, k}$, providing that Conjecture \ref{conj3} holds.
}
\end{remark}

\section{Concluding reflections}\label{conclud}

Study of the interpolation problem for $\hol(\D, \Gamma)$
 was originally motivated by a wish to solve the ``$\mu$-synthesis problem", which arises in control engineering \cite{Do,DuPa}. This is a hard problem of a function-theoretic nature, and its solution would have considerable significance for engineers. Unfortunately, at present it can be analysed in only a few very special cases \cite{NJY11}; in this paper we throw some light on a further case -- the {\em spectral Nevanlinna-Pick problem for $2 \times 2$ matrix functions} with $n > 2$ interpolation points. Given points $\la_1, \dots, \la_n \in \D$ and target matrices $W_1, \dots, W_n \in \C^{2 \times 2}$ one seeks an analytic $2 \times 2$-matrix-valued function $F$ such that 
\[
F(\la_j)= W_j \quad \mbox{  for } j=1,\dots,n,  \mbox{  and }
\]
\[
r(F(\la)) \leq 1 \quad \mbox{ for all } \la \in\D
\]
where $r$ denotes the spectral radius. This problem is essentially equivalent to the interpolation problem for $\hol(\D, \Gamma)$ studied here; see \cite[Theorem 1.1]{AY04T}.

If  Conjecture \ref{conj3} is true then one can check\footnote{subject to a minor complication in the case that some $W_j$ is a scalar matrix.} whether a given spectral Nevanlinna-Pick problem $\la_j \mapsto W_j$, $1 \le j \le n$, has a solution by determining whether the $\Gamma$-interpolation data $\la_j \mapsto ({\rm tr~}  W_j, \det  W_j)$, $1 \le j \le n$, satisfy condition 
$\mathcal{C}_{n -2}$ (where $ n \ge 2$). To verify condition 
$\mathcal{C}_{n -2}$ one must check for positivity a pencil of 
$n \times n$ matrices indexed by $\B l_{n -2}$, the set of  Blaschke products  of degree at most $n-2$. Now $\B l_{n -2}$
is a compact set of real dimension $2n -3$ in the topology of locally uniform convergence. In cases of engineering interest $n$ is likely to be small, and so there is a fair prospect that condition
$\mathcal{C}_{n -2}$ can be checked efficiently. We have not attempted any numerical studies. Engineers currently use a heuristic algorithm called ``D-K iteration"
\cite[Section 9.3.3]{DuPa}, based on results of Bercovici, Foias and Tannenbaum \cite{BFT90,BFT1}, to attempt to solve $\mu$-synthesis problems, but this algorithm is slow and unreliable. For the $n$-point $\Gamma$-interpolation problem it requires a search over an unbounded, nonconvex set of $6n$ real dimensions. At least for this very special case of $\mu$-synthesis, if  the $\Gamma$-interpolation Conjecture is true then one should be able to improve substantially on current methods.

Finding good algorithms is one goal of our research, but equally important is to develop a satisfactory analytic theory of $\mu$-synthesis problems. For example, in proving the $\Gamma$-interpolation
 Conjecture one might be able to show that a solvable $n$-point $\Gamma$-interpolation problem has a solution that is $\Gamma$-inner of degree at most $2n -2$. Aside from its theoretical interest, such a result could have   practical applications.
A good analytic theory would explain the phenomenon of ill-conditioning which engineers have encountered, and would enable numerical analysts to test their algorithms against a range of examples that are exactly solvable. One could also hope to derive parametrizations of solution sets of a range of $\mu$-synthesis
problems, like those that exist for classical Nevanlinna-Pick problems.

We finish with some questions whose answers would be significant for the understanding of $\mu$-synthesis.  The main question we leave open is of course whether the $\Gamma$-interpolation Conjecture is true, but here are two more.
\begin{question}
Is there a bang-bang theorem for $n$-extremals?
\end{question}
That is, if $h$ is $n$-extremal for $\hol(\Omega_1,\Omega_2)$, then does $h$ necessarily map the topological or distinguished boundary of $\Omega_1$ into the corresponding boundary of $\Omega_2$?  We are particularly interested in the question in the case of $\hol(\D,\G)$.  There is a general bang-bang theorem due to J. W. Helton and R. Howe which can be applied to $\hol(\D,\Omega)$ \cite{HH}, but it assumes that  $\Omega$  has a smooth boundary, and so does not apply to $\G$.
\begin{question}
Is every $n$-extremal in $\hol(\D,\Gamma)$ rational?
\end{question}

JIM AGLER, Department of Mathematics, University of California at San Diego, CA \textup{92103}, USA\\

ZINAIDA A. LYKOVA,
School of Mathematics and Statistics, Newcastle University,
 NE\textup{1} \textup{7}RU, U.K.~~
e-mail\textup{: \texttt{Z.A.Lykova@newcastle.ac.uk}}\\

N. J. YOUNG, School of Mathematics, Leeds University, LS2 9JT, U.K.~~
e-mail\textup{: \texttt{N.J.Young@leeds.ac.uk}}
\end{document}